\documentclass[11pt, reqno ]{amsart}
\usepackage{geometry}                % See geometry.pdf to learn the layout options. There are lots.
\usepackage{graphicx}
\usepackage{amssymb}
\usepackage{epstopdf}
\usepackage{amsmath}
\usepackage{mathrsfs}
  \usepackage{array}
  \usepackage{amsthm}
  \usepackage{enumitem}
\usepackage[all]{xy}
\usepackage{color}
\usepackage{tikz}
\usepackage[colorlinks=true, pdfstartview=FitV,linkcolor=blue,citecolor=blue,urlcolor=blue]{hyperref}

    \advance \topmargin by -\headheight
    \advance \topmargin by -\headsep

    \evensidemargin \oddsidemargin
\newcommand{\subscript}[2]{$#1_#2$}

\numberwithin{equation}{section}

\begin{document}

 \newtheorem{theorem}{Theorem}[section]
 \newtheorem{acknowledgement}[theorem]{Acknowledgement}
  \newtheorem{prop}[theorem]{Proposition}
  \newtheorem{cor}[theorem]{Corollary}
  \newtheorem{lemma}[theorem]{Lemma}
  \newtheorem{defn}[theorem]{Definition}
  \newtheorem{ex}[theorem]{Example}
    \newtheorem{conj}[theorem]{Conjecture}
    \newtheorem{rk}[theorem]{Remark}

\title{Chiral differential operators on the upper half plane and modular forms}

\author[Xuanzhong Dai]{Xuanzhong Dai}
\address{Department of Mathematics, The Hong Kong University of Science and Technology, Clear Water Bay, Kowloon, Hong Kong}
\email{ xdaiac@connect.ust.hk    }

\thanks{ }

\maketitle

\maketitle

\section { Introduction }

Consider the Heisenberg Lie algebra with basis $a_n,b_n\;(n\in \mathbb Z )$, the central element $C$, and with commutation relations
\begin{equation}\label{hei1}
[a_m,b_n]=\delta_{m,-n}C.
\end{equation}
Its vacuum representation $V=\mathbb C[a_{-1},a_{-2},\cdots,b_0,b_{-1},\cdots]$ generated by the vacuum vector $1$, with the relations
\[
a_m 1=0 \;\;\;\text{ if } m\geq 0;\;\;\;\; b_n 1=0\;\;\; \text{ if } n>0; \;\;\;\; C1=1,
\]
 has a structure of vertex operator algebra. Let 
\[
\mathbb H:=\{ \tau \in\mathbb C \; | \text{ im }\tau >0\}
\]
be the upper half plane. By the result of Malikov, Schectman and Vaintrob \cite{MSV},
\begin{equation}\label{obj1}
\mathscr D^{\text{ch}}(\mathbb H):=V \otimes_{\mathbb C [b_0]} \mathcal O(\mathbb H),
\end{equation}
where $\mathbb C[b_0]$ is considered as a subring of the ring of holomorphic functions $\mathcal O(\mathbb H)$ on $\mathbb H$ by $b_0 \mapsto \tau$,  is also a vertex operator algebra, which is called the vertex algebra of chiral differential operators on $\mathbb H$.

It can be proved that the $SL(2,\mathbb R)$-action on $\mathbb H$ by the fractional linear transformation induces an action of $SL(2,\mathbb R) $ on $\mathscr D^{\text{ch}}(\mathbb H)$ as automorphisms of vertex algebras (see Section 2). Let $\Gamma(1):=SL(2,\mathbb Z)$, and
$\Gamma\subset \Gamma(1)$ be an arbitrary congruence subgroup. 
In this work, we will study the fixed point vertex algebra $\mathscr D^{\text{ch}}(\mathbb H)^\Gamma$ under $\Gamma$-action. As in the theory of modular forms, we consider the subspace of $\mathscr D^{\text{ch}}(\mathbb H)^\Gamma$ consisting of elements that are holomorphic at the cusps (see Section 2 for definition), denoted by $\mathscr D^{\text{ch}}(\mathbb H,\Gamma)$.
Since the $SL(2,\mathbb R)$-action preserves the conformal weights of $\mathscr D^{\text{ch}}(\mathbb H)$, $\mathscr D^{\text{ch}}(\mathbb H,\Gamma)$ is naturally a $\mathbb Z_{\geq 0}$ graded vertex operator algebra. 
One of the main purposes of this work is to understand the structure of $\mathscr D^{\text{ch}}(\mathbb H,\Gamma)$ and compute its character.

We will show that the structure of $\mathscr D^{\text{ch}}(\mathbb H,\Gamma)$ is closely related to the modular forms of level $\Gamma$.  Let $M_k(\Gamma)$ be the space of modular forms of weight $k$. 
For any $f\in M_{2k}(\Gamma), k >0$, we introduce a certain subspace $D_f$ of $\mathscr D^{\text{ch}}(\mathbb H,\Gamma)$ which is obtained by applying invariant vertex operators to $f$ (see Section 5 for precise definition).
For $f\equiv 1$, $D_1$ is obtained by invariant vertex operators and a quasi-modular form $E_2$ (see also in Section 5).

Let $ \mathscr B$ be a homogeneous basis of $\oplus_{k\geq 0} M_{2k}(\Gamma)$, then we will prove that $\mathscr D^{\text{ch}}(\mathbb H,\Gamma)$ can be decomposed as a direct sum of $D_f$ for $f\in \mathscr B$, i.e.
\[
\mathscr D^{\text{ch}}(\mathbb H,\Gamma)=\oplus_{f\in \mathscr B} D_f.
\]
One of the main results of this work is the following.
\begin{theorem} \label{theorem1.1}
The character formula of $ \mathscr D^{\text{ch}}(\mathbb H,\Gamma)$ is given by 
\[
\sum_{m=0}^\infty \sum_{n= 0}^\infty \dim M_{2m}(\Gamma) q^{2n+m}\prod_{i=1}^n \dfrac{1}{1-q^i} \prod_{j=1}^{m+n}\dfrac{1}{1-q^j}.
\]
\end{theorem}

For any partition $\lambda=(\lambda_{(1)},\cdots,\lambda_{(d)})$, we define $|\lambda|:=\sum_{i=1}^d \lambda_{i}$, $p(\lambda):=d$.
Since $SL(2,\mathbb R)$ preserves the conformal weight, we consider the conformal weight $N$ subspace $\mathscr D^{\text{ch}}(\mathbb H)_N$ of $\mathscr D^{\text{ch}}(\mathbb H)$. 
For each partition pair $(\lambda,\mu)$ such that $|\lambda|+|\mu| =N,$ we will introduce an $SL(2,\mathbb R)$-invariant subspace $ V_{\lambda,\mu}\subset \mathscr D^{\text{ch}}(\mathbb H)_N$. 
And we introduce a total order on the partition pairs such that whenever $(\lambda',\mu')<(\lambda,\mu)$ and $|\lambda'|+|\mu'|=|\lambda|+|\mu|$, we have $V_{\lambda',\mu'}\subset V_{\lambda,\mu} $, 
thus we obtain a filtration labeled by partition pairs in $\mathscr D^{\text{ch}}(\mathbb H)_N $.
Then we will consider a subspace $(V_{\lambda,\mu})^\Gamma_0 \subset V_{\lambda,\mu}$, consisting of $\Gamma$-invariant elements that satisfy the cuspidal conditions (see Section 2). 
Let $(\lambda_1,\mu_1)$ be the largest partition pair under the condition that $(\lambda_1,\mu_1)<(\lambda,\mu)$ and $|\lambda_1|+|\mu_1|=|\lambda|+|\mu|$, 
then another main result is (Theorem \ref{theorem2.5})
\[
(V_{\lambda,\mu})^\Gamma_0 /(V_{\lambda',\mu'})^\Gamma_0 \cong M_{2k}(\Gamma),
\]
where $k=p(\mu)-p(\lambda)$.

The structure of the paper is as follows: 
in Section 2, we will introduce an $SL(2,\mathbb R)$-invariant filtration labeled by partition pairs in $\mathscr D^{\text{ch}}(\mathbb H)_N$ and state Theorem \ref{theorem2.5}.  
In Section 3, we give a proof of Theorem \ref{theorem2.5} for the case $p(\lambda)-p(\mu)\leq -1$.
In Section 4, we give the proof of Theorem \ref{theorem2.5} for the case $p(\lambda)=p(\mu)$, and in this case $M_0(\Gamma)$ only consists of the constant modular forms.
In Section 5, we will describe the structure of $\mathscr D^{\text{ch}}(\mathbb H,\Gamma)$ and compute the character formula. 
In Section 6, we will give an explicit formula
 for the lifting of arbitrary $f\in M_{2k}(\Gamma(1))$ with $k\geq 1$ to $(V_{\lambda,\mu})^{\Gamma(1)}_0$ for the case $\lambda=\emptyset$ and $p(\mu)=k$.

\noindent \textbf{Acknowledgement.}
The author wishes to thank his advisor, Prof. Y. Zhu for discussion.

\section{The algebra of chiral differential operators on the upper half plane}
In this section, we recall the construction of the vertex algebra $\mathscr D^{\text{ch}}(\mathbb H)$ of chiral differential operators on $\mathbb H$ and construct an $SL(2,\mathbb R)$-action as in \cite{MSV}.
And we introduce an $SL(2,\mathbb R)$-invariant filtration on $\mathscr D^{\text{ch}}(\mathbb H)$, and a cuspidal condition on the $\Gamma$-fixed algebra $\mathscr D^{\text{ch}}(\mathbb H)^{\Gamma}$.

The vacuum representation $V$ in Section 1 of the Heisenberg Lie algebra (\ref{hei1}) is a polynomial algebra of variables $b_0,b_{-1},\cdots,a_{-1},a_{-2},\cdots$, and the Virasoro element is given by
\[
\omega=a_{-1}b_{-1}.
\]
Then $L_0=\omega_{(1)}$ gives $V$ a gradation $V=\oplus _{n=0}^\infty V_n$, where an element in $V_n$ is said to have conformal weight $n$. So $V_0=\mathbb C[b_0]$ and $V_1=\mathbb C[b_0]a_{-1}\oplus \mathbb C[b_0] b_{-1}$.  We will write $a=a_{-1}\cdot 1,\; b=b_0$.
The basic fields $a(z)$ and $b(z)$ are given by
\[
a(z)= \sum_{n\in \mathbb Z} a_n z^{-n-1},\;\;\;
b(z)=\sum_{n\in \mathbb Z} b_nz^{-n}.
\]
According to \cite{MSV}, $\mathscr D^{\text{ch}}(\mathbb H)$ as in (\ref{obj1}) is also a vertex operator algebra generated by the basic fields, $a(z),b(z)$ as above and $Y(f,z), f\in \mathcal O(\mathbb H)$, where
\begin{equation} \label{2.1}
Y(f,z)=\sum_{i= 0}^\infty \frac{\partial^i}{i!}  f(b) (\sum_{n\neq 0} b_nz^{-n})^i.
\end{equation}
We write $f(b)_{m+1}:=f(b)_{(m)}$ for the coefficient of $z^{-m-1}$ in the field $Y(f,z)$.

Certain vertex operators on $\mathscr D^{\text{ch}}(\mathbb H)$ generates representations of affine Kac-Moody algebra $\widehat{\mathfrak{sl}}_2$. More precisely let
\begin{equation}\label{2.2}
E:=-a_{-1},\; F:= a_{-1}b_0^2 +2b_{-1},\; H:= -2a_{-1}b_0.
\end{equation}
We have the following theorem
\begin{theorem} \cite{W,FF,F} \label{theorem2.1}
The coefficients of $E_{(n)}, F_{(n)}, H_{(n)}$ of fields $Y(E,z), Y(F,z), Y(H,z)$ satisfy the relations of affine Kac-Moody algebra $\widehat{\mathfrak{sl}}_2$ of critical level $-2$, where $E,F,H$ corresponds to matrices 
\[\begin{pmatrix}  0 & 1 \\ 0 & 0\end{pmatrix},\;
 \begin{pmatrix} 0 & 0\\ 1 & 0 \end{pmatrix},\;
  \begin{pmatrix} 1 & 0 \\ 0 & -1 \end{pmatrix}
  \] respectively.
\end{theorem}

This representation of $\widehat{\mathfrak{sl}}_2$ on $V$ was first introduced by M. Wakimoto \cite{W}. The general construction of Wakimoto modules was given by B. Feigin and E. Frenkel \cite{FF}.
Its connection with vertex algebras as in the above formulation can be found in \cite{F}.

It is well-known that there is a natural right $SL(2,\mathbb R)$-action on the space 
\[
\Omega(\mathbb H)=\Omega^0(\mathbb H) \oplus \Omega^1(\mathbb H)=\{ f(b)+g(b)db\}
\]
induced by the fractional linear transformation on $\mathbb H$:
\[
\begin{pmatrix}
\alpha & \beta\\
\gamma & \delta
\end{pmatrix}
\tau=
\dfrac{\alpha \tau+\beta}{\gamma \tau+\delta}.
\]
So the Lie algebra $\mathfrak{sl}_2$ acts on $\Omega(\mathbb H)$ as Lie derivatives, where 
$\begin{pmatrix}  0 & 1 \\ 0 & 0\end{pmatrix}, \begin{pmatrix} 0 & 0\\ 1 & 0 \end{pmatrix},
  \begin{pmatrix} 1 & 0 \\ 0 & -1 \end{pmatrix}$ act as $-\frac{d}{db}, b^2\frac{d}{db}, -2b \frac{d}{db}$ respectively.
  And we replace $\frac{d}{db}$ by $a_{-1}$, we get the formula of $E$ and $H$ in (\ref{2.2}). But for the formula of $F$, we need to add an extra term $2b_{-1}$.

For a vector $v$ of conformal weight $1$ in a vertex algebra with field $Y(v,z)=\sum_{n\in \mathbb Z} v_{(n)}z^{-n-1}$, $v_{(0)}$ is a derivation.
In our case $E_{(0)},F_{(0)}$ and $H_{(0)}$ give an action of Lie algebra $\mathfrak{sl}_2$ on $\mathscr D^{\text{ch}}(\mathbb H)$ as derivations. 
By the method in \cite{MSV}, we can show that this can be integrated to an $SL(2,\mathbb R)$-action as automorphisms of vertex algebra.
Because we will consider the action of a congruence subgroup $\Gamma \subset SL(2,\mathbb R)$, and it will be related to the theory of modular forms, where the action of $SL(2,\mathbb R)$ is always from the right (see, e.g., \cite{B}), 
we will make our action of $SL(2,\mathbb R)$ a right action. 
By definition, for $g=e^{x}\in SL(2,\mathbb R),\; x\in \mathfrak{sl}_2$, then
\[
\pi(g)=\sum_{n\geq 0} \frac{(-x_{(0)})^n}{n!}.
\]
And we have $\pi(g_1g_2)=\pi(g_2)\pi(g_1)$.

The $SL(2,\mathbb R)$-action commutes with the translation operator $T=L_{-1}=\omega_{(0)}$ for the fact that
\[
[T,x_{(0)}]=(Tx)_{(0)}=0,\;\;\;\text{ for } x\in\mathfrak{sl}_2\subset \mathscr D^{\text{ch}}(\mathbb H).
\]
 And it also commutes with the semisimple operator $L_0=\omega_{(1)}$, so it preserves the gradation.

We now give the formula of the action of 
 \begin{equation} \label{2.3}
 g=\begin{pmatrix}
 \alpha  & \beta\\
 \gamma &\delta
 \end{pmatrix}\in SL(2,\mathbb R)
 \end{equation} 
on generators $a,b,f(b)\in \mathscr{D}^{ch}(\mathbb H)$

\begin{align} 
\nonumber \pi(g) a&=a_{-1}(\gamma b+\delta)^2 +2 \gamma^2 b_{-1}\\
\label{2.4}\pi(g) b&=\dfrac{\alpha b+\beta}{\gamma b +\delta}\\
\nonumber\pi(g)f(b)&=f(gb)=f(\dfrac{\alpha b+\beta}{\gamma b +\delta})
\end{align}

For simplicity, we will introduce a notation $a_{-\lambda}$ for a long expression 
\[
a_{-\lambda_{(1)}}a_{-\lambda_{(2)}}\cdots a_{-\lambda_{(d)}}, 
\]
where $\lambda=(\lambda_{(1)},\cdots, \lambda_{(d)})$ with $ \lambda_{(1)}\geq \lambda_{(2)}\geq \cdots \geq \lambda_{(d)}\geq 1$, i.e. $\lambda$ is a partition.  And we define $p(\lambda):=d$, and $|\lambda |:=\sum_{i=1}^d \lambda_{(i)}$. 
So $\lambda$ is a partition of $|\lambda|$ with $p(\lambda)$ parts. 
Similarly we define $b_{-\mu}$ for any partition $\mu$.
Every element in $\mathscr D^{\text{ch}}(\mathbb H)$ can be written as a sum of elements of type $a_{-\lambda}b_{-\mu}f(b)$ with $f(b)\in \mathcal O(\mathbb H)$. Notice that $a_{-\lambda}b_{-\mu} f(b)$ has conformal weight $|\lambda|+|\mu|$.
We also consider the empty set as a partition, and set $a_{-\emptyset} =b_{-\emptyset}=1$, and $|\emptyset|=p(\emptyset)=0$.

Notice that as a subgroup of $SL(2,\mathbb R)$, the congruence subgroup $\Gamma\subset SL(2,\mathbb Z)=\Gamma(1)$ also acts on $\mathscr D^{\text{ch}}(\mathbb H)$. We denote by $\mathscr D^{\text{ch}}(\mathbb H)^\Gamma$ the $\Gamma$-fixed points of $\mathscr D^{\text{ch}}(\mathbb H)$.
$\mathscr{D}^{\text{ch}}(\mathbb H)$ is not an interesting object as it is too big, so we consider the elements in $\mathscr{D}^{\text{ch}}(\mathbb H)^\Gamma$ satisfying the cuspidal conditions similar to the definition of modular forms of $\Gamma$.

We consider $\Gamma=\Gamma(1)$ first. 
Since $\pi\left(\begin{pmatrix} 1 & 1\\0 &1\end{pmatrix} \right)a=a,\; \pi\left(\begin{pmatrix} 1 & 1\\0 &1\end{pmatrix} \right) b=b+1$ by (\ref{2.4}), $\pi \left(\begin{pmatrix} 1 & 1\\0 &1\end{pmatrix} \right) $ preserves $a_{-n}$ and $b_{-n}$ for $n\geq 1$. 
And it acts as an automorphism on $\mathscr D^{\text{ch}}(\mathbb H)$, so we have
\[
\pi \left(\begin{pmatrix} 1 & 1\\0 &1\end{pmatrix} \right) \sum a_{-\lambda}b_{-\mu} f_{\lambda,\mu}(b)=\sum a_{-\lambda}b_{-\mu} f_{\lambda,\mu}(b+1).
\]
Hence $f_{\lambda,\mu}(b+1)=f_{\lambda,\mu}(b)$, and $f_{\lambda,\mu}(b)$ has a $q$-expansion at the cusp $\infty$,
 \[
 f_{\lambda,\mu}(b)= \sum_{m=-\infty}^{\infty} u_{\lambda,\mu}(m) q^m , \;\;\;\text{ where } q=e^{2\pi i b}.
 \]
 We call $v=\sum_{\lambda,\mu} a_{-\lambda} b_{-\mu} f_{\lambda,\mu}$ is holomorphic at $\infty$, if 
  for arbitrary partitions $\lambda,\mu$, we have $u_{\lambda,\mu}(m)=0$ for $m<0$. Since all the cusps $\mathbb Q\cup \{\infty\}$ are $SL(2,\mathbb Z)$-equivalent, we call $v$ is holomorphic at the cusps.
 
 For a general congruence subgroup $\Gamma$, the notion of holomorphicity at the cusp $c\in \mathbb Q\cup \{\infty\}$ needs more discussions. 
 Choose $\rho\in SL(2,\mathbb Z)$ such that $\rho (c)=\infty$.  
 Then $\pi(\rho)v=\sum a_{-\lambda'}b_{-\mu'}\tilde {f}_{\lambda',\mu'}$ is invariant under $\rho^{-1} \Gamma \rho$ as the group action is a right action. 
 And since $\rho^{-1} \Gamma \rho$ contains the translation matrix $\begin{pmatrix}  1 & N \\ 0 & 1 \end{pmatrix}$,  for some positive integer $N$(cf. \cite{B} p.41-42), $\pi(\rho)v$ is fixed by $\begin{pmatrix}  1 & N \\ 0 & 1 \end{pmatrix}$,
 which implies that $\tilde{f}_{\lambda',\mu'}(b_0+N)=\tilde{f}_{\lambda',\mu'}(b_0)$.
 Hence $\tilde f_{\lambda',\mu'}$ has a Fourier expansion $\sum \tilde{u}_{\lambda',\mu'}(m) e^{2\pi i mb/N}$. 
We say that $v$ is holomorphic at the cusp $c$ if for arbitrary partitions $\lambda',\mu'$, we have $\tilde{u}_{\lambda',\mu'}(m)=0$ for $m<0$. 
 We use $\mathscr D^{\text{ch}}(\mathbb H,\Gamma)$ to denote the $\Gamma$-invariant vectors in $\mathscr D^{\text{ch}}(\mathbb H)$ that are holomorphic at all the cusps. Using (\ref{2.1}), we can prove $\mathscr{D}^{\text{ch}}(\mathbb H,\Gamma)$ is a vertex subalgebra. 
 
 \begin{prop}
$\mathscr D^{\text{ch}}(\mathbb H,\Gamma) $ is a vertex subalgebra of $\mathscr D^{\text{ch}}(\mathbb H)$.
\end{prop}

For $x\in \mathscr D^{\text{ch}}(\mathbb H,\Gamma)$, the adjoint action of $g\in SL(2,\mathbb R)$ on the operator $x_{(n)}$ is defined to be
\begin{equation}
\pi(g) x_{(n)} \pi(g)^{-1}=(\pi(g)x)_{(n)}.
\end{equation}
In particular, the formulas of the adjoint action on the operators $a_{-n}$ and $b_{-n}$ for $n\geq 1$ are given by
\begin{align} \nonumber
\pi(g)a_{-n} \pi(g)^{-1}&=(\pi(g)a)_{-n} =(a_{-1} (\gamma b+\delta)^2)_{-n} +2n\gamma^2 b_{-n}\\
\nonumber
&=  \sum_{k\geq 1} a_{-k} (\gamma b+\delta)^2 _{-n+k} + \sum_{k\geq 0} (\gamma b+\delta)^2_{-n-k}a_k  +2n \gamma^2 b_{-n} \\
\nonumber &=a_{-n} \left((\gamma b+\delta)^2+\gamma^2 \sum_{i\neq 0} b_{-i}b_i  \right) +\sum_{\substack{k\geq 1,\\k\neq n}} a_{-k} \left(2\gamma (\gamma b+\delta) b_{-n+k} +\gamma ^2\sum_{\substack{i,j\neq 0\\ i+j=n-k}} b_{-i}b_{-j}  \right)\\
\label{2.6}
& +\sum_{k\geq 0} \left(2\gamma (\gamma b+\delta) b_{-n-k} +\gamma^2 \sum_{\substack{i,j\neq 0\\ i+j=n+k}}b_{-i}b_{-j}  \right)a_k +2n \gamma^2 b_{-n}\\
 \nonumber
 \pi(g)b_{-n} \pi(g)^{-1}&=  (\pi(g)b)_{-n}=(\dfrac{\alpha b+\beta}{\gamma b+\delta})_{-n}\\
 \label{2.7}
&= \sum_{l\geq 1} \sum_{\substack{i_1,\cdots,i_l   \in \mathbb Z_{\neq 0}  :\\ i_1+\cdots+i_l=n }} (-\gamma)^{l-1}(\gamma b+\delta)^{-l-1} b_{-i_1}\cdots b_{-i_l}
 \end{align}
where the third equality in (\ref{2.6}) is given by the Borcherds identity, and the last equalities in (\ref{2.6}) and (\ref{2.7}) are given by (\ref{2.1}).
Since $SL(2,\mathbb R)$ acts on $\mathscr D^{\text{ch}}(\mathbb H)$ as automorphisms, the action of $g$ is given by
\begin{equation} \label{2.8}
\pi(g) a_{-\lambda} b_{-\mu} f(b) =(\pi(g) a)_{-\lambda_{(1)}}\cdots (\pi(g)a)_{-\lambda_{(n)}} (\pi(g)b)_{-\mu_{(1)}}\cdots (\pi(g)b)_{-\mu_{(m)}}f(gb),
\end{equation}
where $\lambda=(\lambda_{(1)},\cdots,\lambda_{(n)}) ,\mu=(\mu_{(1)},\cdots,\mu_{(m)})$ are partitions.

If we replace $(\pi(g)a)_{-\lambda_{(i)}}$ and $(\pi(g) b)_{-\mu_{(j)}}$ in the right side of (\ref{2.8}) by (\ref{2.6}) and (\ref{2.7}), we can prove that
 \begin{lemma}\label{lemma2.3}
 For $g$ as in (\ref{2.3}), and holomorphic function $f$ on $\mathbb H$, 
\begin{equation}\label{2.9}
\pi(g) a_{-\lambda} b_{-\mu} f(b) =a_{-\lambda}b_{-\mu} (\gamma b+\delta)^{2(p(\lambda)-p(\mu))}f(gb)+\sum_{\substack{\lambda',\mu':  p(\lambda')\leq p(\lambda)  \\p(\lambda')-p(\mu')<p(\lambda)-p(\mu),\\ |\lambda'|+|\mu'|=|\lambda|+|\mu|}} a_{-\lambda'}b_{-\mu'} f_{\lambda',\mu'}(b),
\end{equation}
where $f_{\lambda',\mu'}$ is a holomorphic function on $\mathbb H$.
\end{lemma}

\noindent {\it Proof: }
After replacing $(\pi(g)a)_{-\lambda_{(i)}}$ and $(\pi(g) b)_{-\mu_{(j)}}$ in the right side of (\ref{2.8}) by (\ref{2.6}) and (\ref{2.7}), and moving the annihilation operators to the right, the result is a sum of elements of type 
 $a_{-\lambda'}b_{-\mu'}f_{\lambda',\mu'}$. If $a_{-\lambda'}b_{-\mu'}f_{\lambda',\mu'}$ appears,
 then $|\lambda'|+|\mu'|=|\lambda|+|\mu|$, because the action preserves the conformal weight. 
 
 Since the formula of the adjoint action on $a_{-n}$ in (\ref{2.6}) has at most one $a_m (m\in \mathbb Z)$ and the adjoint action on $b_{-n}$ in (\ref{2.7}) is free of $a_m(m\in \mathbb Z)$, so the $\pi(g)$ action will not increase the number of $a_m(m\in \mathbb Z)$, namely $p(\lambda')\leq p(\lambda)$. 
 Notice that $p(\lambda')-p(\mu')$ means the difference of the number of $a_{-m}(m\geq 1)$ and the number of $b_{-m}(m\geq 1)$ in each expression $a_{-\lambda'}b_{-\mu'}f_{\lambda',\mu'}$. 

From (\ref{2.6}) and (\ref{2.7}), we can see that all the terms except $a_{-n}(\gamma b+\delta)^2$ in (\ref{2.6}) and $(\gamma b+\delta)^{-2} b_{-n}$, the case when $l=1$ in (\ref{2.7}), will decrease the difference of the number of $a_{-m}(m\geq 1)$ and the number of $b_{-m}(m\geq 1)$. For example the term $\gamma ^2 a_{-n}b_{-i}b_i (i\neq 0)$ in (\ref{2.6}) will decrease the number of $a_{-m}(m\geq 1)$ by $1$, and increase the number of $b_{-m}(m\geq 1)$ by $1$;  the term $2\gamma (\gamma b+\delta)b_{-n-k} a_k (k\geq 0)$ will decrease the number of $a_{-m}(m\geq 1)$ by $1$, and it will preserve the number of $b_{-m}(m\geq 1)$ when $k> 0$, and increase the number of $b_{-m}(m\geq 1)$ by $1$ when $k=0$.
Hence we have
 \[ 
p(\lambda')-p(\mu')\leq p(\lambda)- p(\mu),
 \]
where the equality holds only for the case $(\lambda',\mu')=(\lambda,\mu)$ and the corresponding term equals 
 \[
 \pushQED{\qed}
 a_{-\lambda}b_{-\mu} (\gamma b+\delta)^{2(p(\lambda)-p(\mu))} f(gb).\qedhere
 \popQED
 \]

Note that a conceptual explanation of the above lemma would be using the infinitesimal adjoint action on the operator $a_{-\lambda}b_{-\mu}$, which is a maximal vector (killed by $E_{(0)}$) and the eigenvalue of which under the action of semisimple operator $H_{(0)}$ equals $2(p(\lambda)-p(\mu))$. And the action of $F_{(0)}$ will strictly lower the $H_{(0)}$-weight by $\mathfrak{sl}_2$-theory, hence the infinitesimal action will not increase the $H_{(0)}$-weight. The action will be discussed in detail in Section 3.

Now we will introduce a total order on the partitions and partition pairs to equip a filtration on $\mathscr{D}^{\text{ch}}(\mathbb H)$.

For two partitions $\lambda,\lambda'$,  we say $\lambda>\lambda'$ if either $\lambda_{(i)}=\lambda'_{(i)}$ for $1\leq i\leq j-1$ and $\lambda_{(j)}>\lambda'_{(j)}$,  where $j\leq \min\{p(\lambda),p(\lambda')  \}$;
or $p(\lambda)>p(\lambda')$ and $\lambda_{(i)}=\lambda'_{(i)}$ for $1\leq i \leq p(\lambda')$.
Hence $\emptyset$ is strictly less than any partitions except itself.

And we say $(\lambda,\mu) > (\lambda', \mu')$, 
if one of the following conditions holds
\begin{enumerate} [label=(\subscript{A}{{\arabic*}})]
\item $p(\lambda)-p(\mu)> p(\lambda')-p(\mu')$; \label{A1}\\
\item $p(\lambda)-p(\mu)=p(\lambda')-p(\mu')$, and $ \lambda > \lambda'$;  \label{A2}\\
\item $\lambda=\lambda', p(\mu)=p(\mu')$ and $\mu<\mu'$.  \label{A3}
\end{enumerate}
Obviously this gives a total order on the partition pairs.

Now we define a free $\mathcal O(\mathbb H)$-module of finite rank for a partition pair $(\lambda,\mu)$:
 \[
 V_{\lambda,\mu}:= \text{Span}_\mathbb C \{ a_{-\lambda'}b_{-\mu'}f(b_0)\in \mathscr D^{\text{ch}}(\mathbb H) | (\lambda',\mu') \leq (\lambda,\mu),  |\lambda'|+|\mu'|=|\lambda|+|\mu|  \}.
 \]
Then by Lemma \ref{lemma2.3}, $V_{\lambda,\mu}$ is an $SL(2,\mathbb R)$-submodule of $\mathscr D^{\text{ch}}(\mathbb H)_N$ for $N=|\lambda|+|\mu|$, where $\mathscr D^{\text{ch}}(\mathbb H)_N$ is the conformal weight $N$ subspace of $\mathscr D^{\text{ch}}(\mathbb H)$.
So we have a filtration of submodules $\{V_{\lambda,\mu}\}$, satisfying that
\[
V_{\lambda',\mu'}\subset V_{\lambda,\mu} \;\;\;\text{ if } (\lambda',\mu')<(\lambda,\mu), |\lambda'|+|\mu'|=|\lambda|+|\mu|.
\]
The filtration of chiral differential operators of different types can be found in \cite{MSV}, \cite{S}.

Given partitions $\lambda_0,\mu_0$, there are only finitely many partition pairs $(\lambda,\mu)$ such that $(\lambda,\mu)<(\lambda_0,\mu_0)$ and $|\lambda|+|\mu|=|\lambda_0|+|\mu_0|$. Let $(\lambda_1,\mu_1)$ be the successive partition pair of $(\lambda_0,\mu_0)$ under the above two conditions, namely
\begin{equation} \label{2.10}
(\lambda_1,\mu_1)=\max_{(\lambda,\mu)}\{(\lambda,\mu)<(\lambda_0,\mu_0)\;|\; |\lambda|+|\mu|=|\lambda_0|+|\mu_0|\}.
\end{equation}

Since $V_{\lambda_0,\mu_0}$ and $V_{\lambda_1,\mu_1}$ are preserved under the group action, the quotient space $V_{\lambda_0,\mu_0}/V_{\lambda_1,\mu_1}$ is also an $SL(2,\mathbb R)$-module under the induced group action.
We have an exact sequence of $SL(2,\mathbb R)$-modules
\[
0\longrightarrow V_{\lambda_1,\mu_1} \longrightarrow V_{\lambda_0,\mu_0}\longrightarrow V_{\lambda_0,\mu_0}/ V_{\lambda_1,\mu_1} \longrightarrow 0.
\]
Taking the $\Gamma$-fixed points of the above sequence, we have the exact sequence
\[
0\longrightarrow V_{\lambda_1,\mu_1}^\Gamma \longrightarrow V_{\lambda_0,\mu_0} ^\Gamma \longrightarrow (V_{\lambda_0,\mu_0}/ V_{\lambda_1,\mu_1})^{\Gamma}.
\]

By (\ref{2.9}), all the terms of $\pi(g)  a_{-\lambda_0} b_{-\mu_0} f(b)$ are contained in $V_{\lambda_1,\mu_1}$ except for $a_{-\lambda_0}b_{-\mu_0}(\gamma b+\delta)^{2(p(\lambda_0)-p(\mu_0))} f(gb)$.
Hence
 $a_{-\lambda_0} b_{-\mu_0} f(b) +V_{\lambda_1,\mu_1}$ is fixed by $\Gamma$, if and only if 
\begin{equation} \label{2.11}
f(b)=(\gamma b+\delta)^{2(p(\lambda_0)-p(\mu_0))} f(gb),\;\;\text{ for any } g\in \Gamma.
\end{equation}
 
 Note that (\ref{2.11}) is the main condition for modular forms of level $\Gamma$. We denote by $(V_{\lambda,\mu})^\Gamma_0$ the subspace of $V_{\lambda,\mu}^\Gamma$ consisting of elements holomorphic at all the cusps. 
 And for any partition pair $(\lambda,\mu)$, we define 
 \begin{equation} \label{2.12}
 l(\lambda,\mu):=p(\lambda)-p(\mu).
 \end{equation}
 So we have shown that
\begin{lemma}\label{lemma2.4} Given any partition pair $(\lambda_0,\mu_0)$ and take $(\lambda_1,\mu_1)$ as in (\ref{2.10}), then
\[
(V_{\lambda_0,\mu_0})^\Gamma_0 /(V_{\lambda_1,\mu_1} )^\Gamma_0\subset M_{-2l(\lambda_0,\mu_0)}(\Gamma)=M_{2(p(\mu_0)-p(\lambda_0))}(\Gamma).
\]
\end{lemma}

\noindent Let $\alpha:
 (V_{\lambda_0,\mu_0} )^\Gamma_0  \longrightarrow   M_{-2l(\lambda_0,\mu_0)}(\Gamma)$ be the map defined by 
 \[ 
 \sum_{(\lambda,\mu)\leq (\lambda_0,\mu_0)} a_{-\lambda}b_{-\mu} f_{\lambda,\mu}  \longmapsto   f_{\lambda_0,\mu_0}.
  \]
 
 \noindent Our first main result is
 
  \begin{theorem} \label{theorem2.5}
For any two successive partition pairs $(\lambda_0,\mu_0)>(\lambda_1,\mu_1)$ as in (\ref{2.10}), we have the short exact sequence:
\begin{equation}\label{2.13}
 0\longrightarrow (V_{\lambda_1,\mu_1})^\Gamma_0 \longrightarrow (V_{\lambda_0,\mu_0})^\Gamma_0 \overset{\alpha}{\longrightarrow} M_{-2l(\lambda_0,\mu_0)}(\Gamma)\longrightarrow 0. 
 \end{equation}
 \end{theorem}
The proof of this theorem will be given in Section 3 and Section 4.
  As a direct corollary of Lemma \ref{lemma2.4}, we have
   \begin{prop} 
 For any congruence subgroup $\Gamma$, we have
 \[
 \dim \mathscr D^{\text{ch}}(\mathbb H,\Gamma)_N<\infty, \;\;\text{ for any } N\geq 0,
 \]
 where $\mathscr D^{\text{ch}}(\mathbb H,\Gamma)_N$ denotes the conformal weight $N$ subspace of $\mathscr D^{\text{ch}}(\mathbb H, \Gamma)$.
 \end{prop}
 
\noindent {\it Proof: }
By Lemma \ref{lemma2.4}, the dimension of $(V_{\lambda_0,\mu_0}) ^\Gamma_0$ is bounded by  $\dim (V_{\lambda_1,\mu_1})^\Gamma_0 +\dim M_{-2l(\lambda_0,\mu_0)}(\Gamma)$.
And notice that $(\emptyset,\underbrace{(1,1,\cdots,1)}_N)\leq (\lambda,\mu)$ for any partition pair with $|\lambda|+|\mu|=N$.
We claim that $(V_{\emptyset,\underbrace{(1,1,\cdots,1)}_N})^\Gamma_0$ is finite dimensional. Indeed, by (\ref{2.4}) and (\ref{2.8}), we have
\[
\pi(g)  b_{-1}^N f(b)=b_{-1}^N (\gamma b +\delta)^{-2N} f(gb),
\]
then $b_{-1}^N f(b)$ is fixed by $\Gamma$ and holomorphic at the cusps if and only if $f$ is a modular form of weight $2N$. 
Thus 
\[
(V_{\emptyset,\underbrace{(1,1,\cdots,1)}_N})^\Gamma_0\cong M_{2N}(\Gamma).
\]

Hence we prove the result by induction.\qed

\

\

\section{Lifting of nonconstant modular forms}

 In this section, we will study the lifting under the map $\alpha$ in (\ref{2.13}) of modular forms of positive even weight, to
   $\mathscr D^{\text{ch}}(\mathbb H,\Gamma)$. We will prove Theorem \ref{theorem2.5} when $l(\lambda,\mu) \leq -1$.
   
 Let $\mathcal U$ be the quotient of universal enveloping algebra of the Heisenberg Lie algebra (\ref{hei1}) by the ideal generated by $C-1$, so $\mathcal U$ is a graded algebra with the gradation given by the conformal weight.
Define a topology on $\mathcal U$ in which a fundamental system of neighborhoods of $0$ consists of the left ideals $\mathcal U_n $ generated by the elements with conformal weight less or equal to $-n$. 
Then $\{\mathcal U_n\}_{n=0}^\infty$ is a decreasing series with the condition that $\cap_{n\geq 0} \mathcal U_n=\{0\}$.
Let $\bar{\mathcal U}$ denote the completion of $\mathcal U$ with respect to the topology (see similar constructions in \cite{FZ}). 
And $\bar{\mathcal U}$ has a fundamental system $\{\bar{\mathcal U}_n\}_{n=0}^\infty$ of neighborhoods of $0$. 
 Note that $\bar{\mathcal U}$ acts on $\mathscr D^{\text{ch}}(\mathbb H,\Gamma)$.

For any $N\geq 0$, we can check that all but finitely many terms in (\ref{2.6}) and (\ref{2.7}) are contained in $\bar{\mathcal U}_N$. 
Thus the adjoint action of $g$ on $a_n$ and $b_n$ are contained in $\bar{\mathcal U}$ and hence the Lie group $SL(2,\mathbb R)$ acts on $\bar{\mathcal U}$.  Let $K$ be the left ideal in $\bar{\mathcal U}$ generated by elements $a_{n}$ and $b_{n}$ for $n\geq 1$. 
Then $K$ is preserved under the $SL(2,\mathbb R)$-action by (\ref{2.6}) and (\ref{2.7}). 
Therefore $\bar{\mathcal U}/K$ has an $SL(2,\mathbb R)$-module structure. For $f\in \mathcal O(\mathbb H) \subset \mathscr D^{\text{ch}}(\mathbb H)$, because $a_nf=b_nf=0$ for $n\geq 1$, so $Kf =0$. 
 Therefore we have a map 
 \begin{align*}
 \bar{\mathcal U}/K \times \mathcal O(\mathbb H) & \longrightarrow \mathscr D^{\text{ch}}(\mathbb H,\Gamma)\\
 (u+K) f&\longmapsto uf,
 \end{align*}
 which is $SL(2,\mathbb R)$-equivariant in the sense that 
\[
\pi(g) Af =(\pi(g) A \pi(g)^{-1}) \pi(g) f, \;\;\;\text{ for any } A\in \bar{\mathcal U}/K,  g\in SL(2,\mathbb R). 
\]

 According to PBW theorem, we may write a basis of $\mathcal U$ as 
 \[
 a_{-\lambda}b_{-\mu}a_0^kb_0^l a_{\lambda'}b_{\mu'},\;\;\text{ for all partitions } \lambda,\mu,\lambda',\mu', \text{ and } k,l\in \mathbb Z_{\geq 0},
 \]
 where $a_{\lambda'}$ denotes the expression $a_{\lambda'_{(1)}}\cdots a_{\lambda'_{(d)}}$, with $d=p(\lambda')$, and similarly for $b_{\mu'}$.  Hence $\bar{\mathcal U}/K$ has the following basis
 \[
S:=\{ a_{-\lambda}b_{-\mu}a_0^kb_0^l \;|\; \text{ for all partitions } \lambda,\mu, \text{ and } k,l\in \mathbb Z_{\geq 0}\}.
 \]

Now we will give an order on $S$. We say $a_{-\lambda_1} b_{-\mu_1} a_0^{k_1} b_0^{l_1} > a_{-\lambda_2} b_{-\mu_2} a_0^{k_2} b_0^{l_2}$ if one of the following conditions holds 
 \begin{enumerate}  [label=(\subscript{B}{{\arabic*}})]
 \item $(\lambda_1,\mu_1)>(\lambda_2,\mu_2)$, where the order is defined as in \ref{A1}-\ref{A3} in Section 2;    \label{B1}\\
 \item $(\lambda_1,\mu_1)=(\lambda_2,\mu_2), k_1>k_2$;   \label{B2}\\
 \item $(\lambda_1,\mu_1)=(\lambda_2,\mu_2), k_1=k_2, l_1<l_2.    \label{B3}$
 \end{enumerate}

 We will first study the lifting of a nonconstant modular form $f$ of weight $-2l(\lambda_0,\mu_0)$ with $l(\lambda_0,\mu_0) \leq -1$ in (\ref{2.13}) to $(V_{\lambda_0,\mu_0})^\Gamma_0$, 
the idea is to find an operator
\begin{equation} \label{3.1}
A=a_{-\lambda_0}b_{-\mu_0}+l.o.t \in \bar{\mathcal U}/K,
\end{equation}
where $l.o.t$ refers to terms which are strictly less than $a_{-\lambda_0}b_{-\mu_0} $,  such that
 \[
 \pi(g)A f=Af ,\; \;\;\text{ for any } g\in \Gamma.
 \] 
 
 The left side of above equation equals 
 \[
 (\pi(g)A \pi(g)^{-1}) f(gb)=(\pi(g)A\pi(g)^{-1})(\gamma b+\delta)^{-2l(\lambda_0,\mu_0)}f(b).
 \] 
 So it suffices to find solutions of operator $A$ of the form (\ref{3.1}), such that 
 \begin{equation}\label{3.2}
 \pi(g)A\pi(g)^{-1}=A(\gamma b+\delta)^{2l(\lambda_0,\mu_0)}, \;\;\; \text{ for any } g\in  SL(2,\mathbb R).
 \end{equation}
  Considering the infinitesimal action, (\ref{3.2}) is equivalent to the following system
 \begin{align} 
 \label{3.3} E_{(0)}. A&=0,\\
\label{3.4} H_{(0)}.A &=2l(\lambda_0,\mu_0) A,\\
\label{3.5} F_{(0)}.A &=-2l(\lambda_0,\mu_0)A b_0,
 \end{align}
 where for $x\in \mathfrak{sl}_2$, we denote by $x_{(0)}.$ the infinitesimal adjoint action of $x$ on $\bar{\mathcal U}/K$, which is given by:  
 \begin{equation}\label{3.6}
x_{(0)}. B= x_{(0)}B-Bx_{(0)},\;\;\; \text{ for any } B\in \bar{\mathcal U}/K.
 \end{equation}
 
\noindent The formulas of $E_{(0)}$ and $H_{(0)}$ on an operator $a_{-\lambda}b_{-\mu} a_0^k b_0^l$ can be calculated easily, 
\begin{gather}\label{3.7}
E_{(0)}. a_{-\lambda}b_{-\mu} a_0^k b_0^l  =-l a_{-\lambda}b_{-\mu}a_0^k b_0^{l-1}, \\
\label{3.8}H_{(0)}. a_{-\lambda}b_{-\mu} a_0^k b_0^l  =2(p(\lambda)-p(\mu)+k-l) a_{-\lambda}b_{-\mu} a_0^k b_0^l .
\end{gather}

 We will construct an $\mathfrak {sl}_2$-submodule of $\bar{\mathcal U}/K$ for the partition pair $(\lambda_0,\mu_0)$. We first define a subset $\mathcal S_{\lambda_0,\mu_0} \subset S$ as follows:
 \begin{equation} \label{3.9}
 \mathcal S_{\lambda_0,\mu_0}:=\{ a_{-\lambda} b_{-\mu} a_0^k b_0^l \in S  \, | \, a_{-\lambda} b_{-\mu} a_0^k b_0^l \leq a_{-\lambda_0}b_{-\mu_0},  |\lambda|+|\mu|=|\lambda_0|+|\mu_0|,     l(\lambda,\mu)+k \leq l(\lambda_0,\mu_0)\},
 \end{equation}
 where $l(\lambda,\mu)$ and $l(\lambda_0,\mu_0)$ as in (\ref{2.12}).
According to (\ref{3.8}), the last condition $ l(\lambda,\mu)+k \leq l(\lambda_0,\mu_0)$ in (\ref{3.9}) implies the operators in $\mathcal S_{\lambda_0,\mu_0}$ has $H_{(0)}$-eigenvalues less or equal to $2l(\lambda_0,\mu_0)$.
 
  \begin{lemma} \label{lemma3.1}
 If $a_{-\lambda} b_{-\mu} a_0^k b_0^l \in \mathcal S_{\lambda_0,\mu_0}  $ for some $l\geq 0$, then $a_{-\lambda}b_{-\mu} a_0^i b_0^j \in \mathcal S_{\lambda_0,\mu_0}, $ for all $0\leq i\leq k$ and $j\geq 0$.
 \end{lemma}
 \noindent {\it Proof:}
 The case  $(\lambda,\mu) <(\lambda_0,\mu_0)$ is trivial from the definition of $\mathcal S_{\lambda_0,\mu_0}$. 
 For the case $(\lambda,\mu)=(\lambda_0,\mu_0)$, we conclude that $k$ must be zero. 
 And in this case $a_{-\lambda_0}b_{-\mu_0} b_0^j\in \mathcal S_{\lambda_0,\mu_0}$ for all $j\geq 0$.\qed

  We let $\mathcal S_{\lambda_0,\mu_0}^0 \subset \mathcal S_{\lambda_0,\mu_0}$ consisting of elements which contain no $b_0$, for example, if $a_{-\lambda} b_{-\mu} a_0^k   \in \mathcal S_{\lambda_0,\mu_0}$ then it is also contained  in $\mathcal S_{\lambda_0,\mu_0}^0$. 
 By (\ref{3.7}), $E_{(0)}$ kills all elements in $\mathcal S^0_{\lambda_0,\mu_0}$.
  
    Let $M_{\lambda_0,\mu_0}$ be the space spanned by $\mathcal S_{\lambda_0,\mu_0}$, and we view $M_{\lambda_0,\mu_0}$ as a subspace of $\bar{\mathcal U}/K$. Then
 we have the following lemma:
 \begin{lemma}\label{lemma3.2} For any partition pair $(\lambda_0,\mu_0)$, 
 $M_{\lambda_0,\mu_0}$ is stable under the action (\ref{3.6}). And $H_{(0)}$ acts semisimply with maximal weight $2l(\lambda_0,\mu_0)$.
 \end{lemma}
 \noindent {\it Proof:}
 It suffices to show that $M$ is preserved by the operators $E_{(0)},H_{(0)},F_{(0)}\in \mathfrak{sl}_2$.
 According to (\ref{3.7}) and Lemma \ref{lemma3.1}, $M$ is preserved by $E_{(0)}$.
 And $H_{(0)}$ is a semisimple operator by (\ref{3.8}).
 
 To prove $F_{(0)}$ also preserves $M$,
 we consider the action formulas of $F_{(0)}=(a_{-1}b^2+2b_{-1})_{(0)}=(a_{-1}b^2)_{(0)}$ on $a_{-n}$ and $b_{-n}$ for $n\geq 0$. 
 A direct calculation shows $F_{(0)} a=-2 a_{-1}b$, and $F_{(0)} b =b^2$. 
 Then according to Borcherds identity, 
\begin{align} \label{3.10}
F_{(0)}. a_{-n}=(F_{(0)}a)_{-n} &=-2(a_{-1}b^2)_{-n}=-2\sum_{i\geq 0} a_{-1-i} b_{-n+i+1}-2\sum_{i\geq 0} b_{-n-i} a_i \;\;\;\text{ for } n\geq 0, \\ 
\label{3.11}  F_{(0)}.b_{-m}=(F_{(0)} b)_{-m} &=(b^2)_{-m}=2b_0b_{-m}+ \sum_{i,j>0,i+j=m} b_{-i}b_{-j} \mod K  \;\;\text{ for } m\geq 1.
\end{align}
 Using the above formulas, it is not hard to see that the following six types of terms appear in $F_{(0)}. a_{-\lambda}b_{-\mu} a_0^k b_0^l$:
\begin{enumerate}  [label=(\subscript{C}{{\arabic*}})]
 \item $a_{-\lambda'} b_{-\mu'} a_0^k b_0^l $ with $p(\lambda')=p(\lambda), \lambda'<\lambda,p(\mu')=p(\mu)+1, |\lambda'|+|\mu'|=|\lambda|+|\mu|$, \label{C1}\\
 \item $a_{-\lambda'}b_{-\mu'} a_0^k b_0^l$ with $p(\lambda')=p(\lambda)-1,p(\mu')=p(\mu),|\lambda'|+|\mu'|=|\lambda|+|\mu|$, \label{C2}\\
 \item $a_{-\lambda'}b_{-\mu'} a_0^{k+1} b_0^l$ with $p(\lambda')=p(\lambda)-1,p(\mu')=p(\mu)+1,|\lambda'|+|\mu'|=|\lambda|+|\mu|$,  \label{C3}\\
 \item $a_{-\lambda}b_{-\mu'} a_0^k b_0^l$ with $p(\mu')=p(\mu)+1, \mu'<\mu, |\mu'|=|\mu|$, \label{C4}\\
 \item $a_{-\lambda}b_{-\mu} a_0^{k-1} b_0^l$, \label{C5}\\
 \item $a_{-\lambda}b_{-\mu}a_0^k b_0^{l+1}$. \label{C6}
\end{enumerate}
 
 We can easily check that all of above types are contained in $\mathcal S_{\lambda_0,\mu_0}$ if $a_{-\lambda}b_{-\mu}a_0^k b_0^l \in \mathcal S_{\lambda_0,\mu_0}$. Hence $F_{(0)}$ preserves $M_{\lambda_0,\mu_0}$. 
 Therefore $M_{\lambda_0,\mu_0}$ has an $\mathfrak{sl}_2$-module structure. 
\qed
 
 Notice that the vectors in $\mathcal S_{\lambda_0,\mu_0}^0$ form a basis of maximal vectors in $M_{\lambda_0,\mu_0}$.
 And the above proof also shows that $F_{(0)}$ will decrease the order of $a_{-\lambda}b_{-\mu} a_0^k b_0^l$ and hence $\mathfrak {sl}_2$-action will not increase the order of elements in $\mathcal S_{\lambda_0,\mu_0}$, 
 which means every operator $X_{(0)}$ for $X\in\mathfrak{sl}_2$ has lower triangular matrix with respect to decreasing basis in $\mathcal S_{\lambda_0,\mu_0}^0$.
 This gives us a hint to find the solution of (\ref{3.3})-(\ref{3.5}) under the decreasing basis. 
 And since $F_{(0)}$ will decrease the weight by $2$, we will
study the relations of decreasing basis in $\mathcal S_{\lambda_0,\mu_0}^0$ of weight $2l(\lambda_0,\mu_0)$ and $2l(\lambda_0,\mu_0)-2$.
 
 \begin{lemma} \label{lemma3.3}
Let $S_{\lambda_0,\mu_0}^0(m)\subset S_{\lambda_0,\mu_0}^0$ be the weight $m$ subspace of $\mathcal S_{\lambda_0,\mu_0}^0$ under $H_{(0)}$-action, namely
\[
\mathcal S_{\lambda_0,\mu_0}^0(m):= \{a_{-\lambda}b_{-\mu}a_0^k \in \mathcal S_{\lambda_0,\mu_0}^0 \;| \: H_{(0)} a_{-\lambda}b_{-\mu}a_0^k=
m  a_{-\lambda}b_{-\mu}a_0^k\}.
\]
Then we have
\[|\mathcal S_{\lambda_0,\mu_0}^0(2l(\lambda_0,\mu_0))|>|\mathcal S_{\lambda_0,\mu_0}^0(2l(\lambda_0,\mu_0)-2)|.\] 
\end{lemma}
\noindent {\it Proof:}
For any term $a_{-\lambda}b_{-\mu}a_0^k \in \mathcal S_{\lambda_0,\mu_0}^0(2l(\lambda_0,\mu_0)-2)$, we will show that $a_{-\lambda}b_{-\mu}a_0^{k+1} \in \mathcal S_{\lambda_0,\mu_0}^0(2l(\lambda_0,\mu_0))$.
Indeed, if $(\lambda,\mu)<(\lambda_0,\mu_0)$, then $a_{-\lambda}b_{-\mu}a_0^{k+1} <a_{-\lambda_0}b_{-\mu_0}$ by \ref{B1}. 
And the $H_{(0)}$-eigenvalue of $a_{-\lambda}b_{-\mu}a_0^{k+1}$ is $2l(\lambda_0,\mu_0)$ by (\ref{3.8}). So $a_{-\lambda}b_{-\mu}a_0^{k+1}$ is contained in $\mathcal S_{\lambda_0,\mu_0}^0(2l(\lambda_0,\mu_0))$.
 If $(\lambda,\mu)=(\lambda_0,\mu_0)$, the condition $l(\lambda,\mu)+k \leq  l(\lambda_0,\mu_0)$ implies that $k$ must be zero. But $a_{-\lambda}b_{-\mu}=a_{-\lambda_0}b_{-\mu_0}$ has weight $2l(\lambda_0,\mu_0)$, which contradicts to the assumption $a_{-\lambda}b_{-\mu}\in \mathcal S_{\lambda_0,\mu_0}^0(2l(\lambda_0,\mu_0)-2)$.
So we have a one to one correspondence from $\mathcal S_{\lambda_0,\mu_0}^0(2l(\lambda_0,\mu_0)-2)$ to $\mathcal S_{\lambda_0,\mu_0}^0(2l(\lambda_0,\mu_0))$. And notice that $a_{-\lambda_0}b_{-\mu_0}\in \mathcal S_{\lambda_0,\mu_0}^0(2l(\lambda_0,\mu_0))$ is not in the image of the above map. Thus  $|\mathcal S_{\lambda_0,\mu_0}^0(2l(\lambda_0,\mu_0))|>|\mathcal S_{\lambda_0,\mu_0}^0(2l(\lambda_0,\mu_0)-2)|$.
\qed

\begin{prop} \label{prop3.4}
Let $M_{\lambda_0,\mu_0}(m)$ be the weight $m$ subspace of $M_{\lambda_0,\mu_0}$, then
\begin{gather*}
M_{\lambda_0,\mu_0}(2l(\lambda_0,\mu_0))=\text{Span}_{\mathbb C} \mathcal S_{\lambda_0,\mu_0}^0(2l(\lambda_0,\mu_0)),\\
M_{\lambda_0,\mu_0}(2l(\lambda_0,\mu_0)-2)=M_{\lambda_0,\mu_0}(2l(\lambda_0,\mu_0)) b_0 \oplus \text{Span}_{\mathbb C} \mathcal S_{\lambda_0,\mu_0}^0 (2l(\lambda_0,\mu_0)-2).
\end{gather*}
\end{prop}

\noindent {\it Proof:}
Since $2l(\lambda_0,\mu_0)$ is the maximal weight by Lemma \ref{lemma3.2}, $E_{(0)}$ kills $M_{\lambda_0,\mu_0}(2l(\lambda_0,\mu_0))$. So $M_{\lambda_0,\mu_0}(2l(\lambda_0,\mu_0))$ is spanned by the maximal vectors of weight $2l(\lambda_0,\mu_0)$.

For $a_{-\lambda} b_{-\mu} a_0^k b_0^l \in M_{\lambda_0,\mu_0}(2l(\lambda_0,\mu_0)-2)$, we have 
\begin{equation}\label{3.12}
p(\lambda)-p(\mu)+k-l=l(\lambda_0,\mu_0)-1.
\end{equation}
And the last condition of (\ref{3.9}) implies
\begin{equation}\label{3.13}
p(\lambda)-p(\mu)+k \leq p(\lambda_0)-p(\mu_0)=l(\lambda_0,\mu_0).
\end{equation}

From (\ref{3.12}) and (\ref{3.13}), we deduce that $l\leq 1$. If $l=1$, the equality holds in (\ref{3.13}). 
Hence $p(\lambda)-p(\mu)+k=l(\lambda_0,\mu_0)$. By Lemma \ref{lemma3.1}, $a_{-\lambda} b_{-\mu} a_0^k\in M_{\lambda_0,\mu_0}(2l(\lambda_0,\mu_0))$. 
Hence $a_{-\lambda} b_{-\mu} a_0^k b_0 \in M_{\lambda_0,\mu_0}(2l(\lambda_0,\mu_0))b_0$. If $l=0$,   $a_{-\lambda} b_{-\mu} a_0^k$ is already in $\mathcal S_{\lambda_0,\mu_0}^0(2l(\lambda_0,\mu_0)-2)$.
\qed

Consider the composition map
\[
\bar{F}=\pi' \circ F_{(0)}: M_{\lambda_0,\mu_0}(2l(\lambda_0,\mu_0))\longrightarrow M_{\lambda_0,\mu_0}(2l(\lambda_0,\mu_0)-2)
\longrightarrow M_{\lambda_0,\mu_0}(2l(\lambda_0,\mu_0)-2) /M_{\lambda_0,\mu_0}(2l(\lambda_0,\mu_0))b,
\] 
where $F_{(0)}$ stands for the restriction of $F_{(0)}$ to $M_{\lambda_0,\mu_0}(2l(\lambda_0,\mu_0))= \text{Span}_{\mathbb C}\mathcal S_{\lambda_0,\mu_0}^0(2l(\lambda_0,\mu_0))$, and $\pi'$ is the standard quotient map.

\begin{theorem} \label{theorem3.5}
When $l(\lambda_0,\mu_0)=p(\lambda_0)-p(\mu_0)\leq -1$, there exists an element of the type 
\[
a_{-\lambda_0}b_{-\mu_0}+l.o.t
\]
where $l.o.t$ means the terms strictly less than $a_{-\lambda_0}b_{-\mu_0}$, in the kernel of
the map
\[\bar{F}: M_{\lambda_0,\mu_0}(2l(\lambda_0,\mu_0)) \longrightarrow M_{\lambda_0,\mu_0}(2l(\lambda_0,\mu_0)-2)/M_{\lambda_0,\mu_0}(2l(\lambda_0,\mu_0))b.\]
\end{theorem}
\noindent {\it Proof:}
Suppose $a_{-\lambda}b_{-\mu} a_0^k$ is a maximal vector of weight $2l(\lambda_0,\mu_0)$, then $p(\lambda)-p(\mu)+k=l(\lambda_0,\mu_0)$ by (\ref{3.8}). 
Let $a_{-\lambda'}b_{-\mu'}a_0^{k'}$ be another maximal vector of weight $2l(\lambda_0,\mu_0)$ with $k>k'$. Then $p(\lambda)-p(\mu)<p(\lambda')-p(\mu')$ and hence $a_{-\lambda}b_{-\mu}a_0^k< a_{-\lambda'}b_{-\mu'}a_0^{k'}$ by \ref{A1}.

We list all the elements in $\mathcal S^0_{\lambda_0,\mu_0}(2l(\lambda_0,\mu_0))$ in the decreasing order
\begin{align}
\nonumber
a_{-\lambda_0}b_{-\mu_0} > &a_{-\lambda_1}b_{-\mu_1}>\cdots>a_{-\lambda_{l_1}}b_{-\mu_{l_1}}>
a_{-\lambda_{l_1+1}}b_{-\mu_{l_1+1}}a_0>\cdots \\ \label{3.14}
>&a_{-\lambda_{l_i+1}}b_{-\mu_{l_i+1}}a_0^i >\cdots >a_{-\lambda_{l_{i+1}}}b_{-\mu_{l_{i+1}}}a_0^i> \cdots >a_{-\lambda_{l_t}}b_{-\mu_{l_t}}a_0^{t-1},
\end{align}
where $p(\lambda_{l_i+j_i})-p(\mu_{l_i+j_i)}+i=l(\lambda_0,\mu_0)$, for $0\leq i\leq t-1,\; 1\leq  j_i\leq l_{i+1}-l_{i}$, and $l_0=0$.

Then the maximal vectors below
\begin{equation}\label{3.15}
a_{-\lambda_{l_1+1}}b_{-\mu_{l_1+1}}>\cdots 
>a_{-\lambda_{l_i+1}}b_{-\mu_{l_i+1}}a_0^{i-1} >\cdots >a_{-\lambda_{l_{i+1}}}b_{-\mu_{l_{i+1}}}a_0^{i-1}> \cdots >a_{-\lambda_{l_t}}b_{-\mu_{l_t}}a_0^{t-2},
\end{equation}
forms a decreasing basis of $\text{Span}_\mathbb C \mathcal S_{\lambda_0,\mu_0}(2l(\lambda_0,\mu_0)-2)$, which is identified with
the quotient space $M_{\lambda_0,\mu_0}(2l(\lambda_0,\mu_0)-2)/M_{\lambda_0,\mu_0}(2l(\lambda_0,\mu_0))b$ by Proposition \ref{prop3.4}.

Then the representative matrix of $\bar{F}$ under the two bases (\ref{3.14}) and (\ref{3.15}) is of the form
\[
(B , C)
\]
where $B$ is an $(l_t-l_1) \times (l_1+1)$ matrix and $C$ is lower triangular matrix of type $(l_t-l_1)\times (l_t-l_1)$, because by \ref{C1}-\ref{C6}, 
 \begin{equation} \label{3.16}
F_{(0)}. a_{-\lambda_{l_i+j_i}}b_{-\mu_{l_i+j_i}}a_0^i = c_{i,j_i}a_{-\lambda_{l_i+j_i}}b_{-\mu_{l_i+j_i}}a_0^i b_0 +d_{i,j_i}a_{-\lambda_{l_i+j_i}} b_{-\mu_{l_i+j_i}}a_0^{i-1}+l.o.t 
\end{equation}
where $l.o.t$ refers to the terms less than $a_{-\lambda_{l_i+j_i}} b_{-\mu_{l_i+j_i}}a_0^{i-1}$, and $c_{i,j_i},d_{i,j_i}$ are constants. 
And now we will calculate the diagonal on $C$, in other word, 
the coefficients $d_{i,j_i}$ for $1\leq i\leq t-1$, and $1\leq j_i \leq l_{i+1}-l_i$.
Notice that
\begin{equation}\label{3.17}
F_{(0)}. a_{-\lambda_{l_i+j_i}}b_{-\mu_{l_i+j_i}}a_0^i  =(F_{(0)}. a_{-\lambda_{l_i+j_i}})b_{-\mu_{l_i+j_i}}a_0^i  +a_{-\lambda_{l_i+j_i}}(F_{(0)}. b_{-\mu_{l_i+j_i}})a_0^i + a_{-\lambda_{l_i+j_i}}b_{-\mu_{l_i+j_i}}F_{(0)}.a_0^i 
\end{equation}
The terms involving $b_0$ in (\ref{3.10}) and (\ref{3.11}) are  $-2a_{-n}b_0$ and $2b_0b_{-m}$ respectively, and notice that $b_0 a_0^{i}=-i a_0^{i-1}+a_0^ib_0$. 
So the contribution of the coefficients $d_{i,j_i}$ for the first two terms in the right side of (\ref{3.17}) is $2ip(\lambda_{l_i+j_i})-2ip(\mu_{l_i+j_i})$.

And the action of $F_{(0)}$ on the derivation operator $a_0^{i}$ equals
\begin{align}\nonumber
F_{(0)}.a_0^{i} &=-2(\sum_{j=1}^{i}  a_0^{i-j}b_0 a_0^j )\\ \nonumber
&=-2\sum_{j=1}^{i} a_0^{i-j} (-ja_0^{j-1}+a_0^j b_0)\\
\label{3.18}&=i(i+1) a_0^{i-1} -2ia_0^{i} b_0
\end{align}

Hence by (\ref{3.17}) and (\ref{3.18}), 
\[
d_{i,j_i}=2ip(\lambda_{l_i+j_i})-2ip(\mu_{l_i+j_i}) +i(i+1)=i(2l(\lambda_{l_i+j_i},\mu_{l_i+j_i})+i+1).
\]
 This is always nonzero, 
because
 \begin{equation}\label{3.19}
 2l(\lambda_{l_i+j_i},\mu_{l_i+j_i})+i+1\leq 2(l(\lambda_0,\mu_0)-i)+i+1=2l(\lambda_0,\mu_0)-i+1<0,
 \end{equation}
  where the first inequality is because $a_{-\lambda_{l_i+j_i}}b_{-\mu_{l_i+j_i}}a_0^i \in \mathcal S_{\lambda_0,\mu_0}$, and the second inequality is because $l(\lambda_0,\mu_0)\leq -1$.  Thus $C$ is an invertible matrix.

Under the bases (\ref{3.14}) and (\ref{3.15}), consider the linear system $\bar{F}(x_0,x_1,\cdots ,x_{l_t})=0$. Let $X=(x_0,\cdots,x_{l_1})^T \in \mathbb C^{l_1+1}, Y=(x_{l_1+1},\cdots, x_{l_{t}})^T \in \mathbb C^{l_t-l_1}$, then $\bar{F}(x_0,x_1,\cdots ,x_{l_t})= (B,C)(X^T,Y^T)^T=BX+CY$, where $M^T$ refers to the transpose of the matrix $M$. Hence for arbitrary $X\in  \mathbb C^{l_1+1}$, $Y=-C^{-1} BX$ is always solvable.  
We take $x_0=1$, and arbitrary $x_1,\cdots,x_{l_1}$, there exists $x_{l_1+1},\cdots,x_{l_t}$, such that $(x_0,\cdots,x_{l_t})$ is a unique nonzero solution of $\bar{F}(x_0,x_1,\cdots ,x_{l_t})=0$.\qed

Note that the proof of Theorem \ref{theorem3.5} doesn't work for the case $l(\lambda_0,\mu_0)=0$, because in this case the matrix $C$ is not invertible.

\noindent {\it Proof of Theorem \ref{theorem2.5}(For the case $l(\lambda_0,\mu_0)\leq  -1$):}
We take a nonzero solution of $\bar{F}(x_0,x_1,\cdots ,x_{l_t})=0$ with $x_0=1$ and let 
\begin{equation}
A= \sum_{i=0}^{l_t} x_i v_i,
\end{equation}
where we denote by $v_0,\cdots,v_{l_t}$ the decreasing basis in (\ref{3.14}).
And Theorem \ref{theorem3.5} shows that when $l(\lambda_0,\mu_0)\leq -1$, 
\[
F_{(0)}. A\in M_{\lambda_0,\mu_0}(2l(\lambda_0,\mu_0))b.
\]
The coefficient of $a_{-\lambda_{l_i+j_i}}b_{-\mu_{l_i+j_i}}a_0^i b_0$, however only comes from the term $F_{(0)}. a_{-\lambda_{l_i+j_i}}b_{-\mu_{l_i+j_i}}a_0^i$, so it equals the constant $c_{i,j_i}$ in (\ref{3.16}). And $c_{i,j_i}= -2p(\lambda_{l_i+j_i})+2p(\mu_{l_i+j_i})-2i=-2l(\lambda_0,\mu_0)$. 
Therefore we have the required operator  $A$ satisfying (\ref{3.5}). And conditions (\ref{3.3}) and (\ref{3.4}) are automatically satisfied. 

Applying $A$ to $f$, we get the lifting of the modular form $f$ in $(V_{\lambda_0,\mu_0})^\Gamma_0$ with the leading term $a_{-\lambda_0}b_{-\mu_0}f$.  \qed

\
\

\section{Lifting of the constant modular forms}

In this section, we will consider 
 the lifting under $\alpha$ in (\ref{2.13}) of the constant modular forms, namely the case when $l(\lambda_0,\mu_0)=0$.  We will prove Theorem \ref{theorem2.5} when $l(\lambda_0,\mu_0)=0$.
 
 To study the lifting of the constant function, we need
 
\begin{lemma}\label{lemma4.1}
 For $g$ as in (\ref{2.3}), the adjoint action of $g$ on an operator $a_{-\lambda}b_{-\mu}\in \bar{\mathcal U}/K$ is given by
 \begin{equation} \label{4.1}
\pi(g) a_{-\lambda} b_{-\mu} \pi(g^{-1}) =\sum_{(\lambda',\mu')\leq (\lambda,\mu)}  \sum_{s=0}^{[\frac{l(\lambda,\mu)-l(\lambda',\mu')}{2}]}  c_{\lambda,\mu,\lambda',\mu'}^s  \gamma^{l(\lambda,\mu)-l(\lambda',\mu')-s} a_{-\lambda'} b_{-\mu'} (\gamma b+\delta)^{l(\lambda,\mu)+l(\lambda',\mu')+s} a_0^s,
\end{equation}
 where $l(\lambda,\mu),l(\lambda',\mu')$ are as in (\ref{2.12}), $[m]$ denotes the greatest integer less than or equal to $m$, and $c_{\lambda,\mu ,\lambda',\mu'}^{s}$ is a constant
 independent of $g$. Moreover $c_{\lambda,\mu,\lambda,\mu}^0=1$,
  and $c_{\lambda,\mu,\lambda',\mu'}^s=0$ if either $|\lambda'|+|\mu'|\neq |\lambda|+|\mu|$ or $(\lambda',\mu')\neq (\lambda,\mu), l(\lambda',\mu')=l(\lambda,\mu)$.
 \end{lemma}
  
  We give an example of the result. 
 \begin{align*}
 \pi(g) a_{-1}b_{-1} \pi(g^{-1})=& a_{-1}b_{-1}+2\gamma (\gamma b+\delta)^{-1} b_{-2} \\
 & - \gamma^2(\gamma b+\delta)^{-2} b_{-1}^2 \\
 & +2\gamma (\gamma b +\delta)^{-1} b_{-1}^2 a_0
 \end{align*}
 where in this case $(\lambda,\mu)=((1),(1))$, and $c_{\lambda,\mu,\lambda',\mu'}^s$ is nonzero only for the cases:
 when $(\lambda',\mu')=((1),(1)), (\emptyset, (2)), (\emptyset, (1,1))$ and $s=0$; when $(\lambda',\mu')=(\emptyset, (1,1))$ and $s=1$.
 
\noindent {\it Proof of Lemma \ref{lemma4.1}:}
Let $g=\begin{pmatrix} \alpha & \beta \\ \gamma & \delta\end{pmatrix}$, and view $\alpha,\beta,\gamma,\delta$ as symbols. 
It is easy to see 
\[
\pi(g) a_{-\lambda}b_{-\mu} \pi(g)^{-1} \in \mathbb C[(\gamma b+\delta),(\gamma b+\delta)^{-1},  \gamma  ] \otimes _{\mathbb C} \bar{\mathcal U}.
\]
We will do the symbolic computation below. We assign an additive degree (denoted by $h$) on the monomials in $C[ (\gamma b+\delta), (\gamma b+\delta)^{-1},  \gamma  ] \otimes _{\mathbb C} \bar{\mathcal U}$ given by:
\[
h(a_n)=-1,\,h(b_n)=1, h(\gamma^{n})=-n ,h((\gamma b+\delta)^n)=0,h(c)=0, \, \text{ for } n\in \mathbb Z,
\]
where $c$ represents the operator multiplied by a constant $c$,  and the ``additive" means for any monomials $X,Y\in  C[(\gamma b+\delta),(\gamma b+\delta)^{-1},  \gamma  ] \otimes _{\mathbb C} \bar{\mathcal U}$
\begin{equation*}
h( XY)=h(X)+ h(Y).
\end{equation*}

We also call a vector homogeneous if the $h$-degrees of its monomials are the same, and we may enlarge the definition of the $h$-degree to the homogeneous vectors in $C[(\gamma b+\delta),(\gamma b+\delta)^{-1},  \gamma  ] \otimes _{\mathbb C} \bar{\mathcal U}$.

Obviously the commutation of the above operators will not change the $h$-degree. For example, the $h$-degrees of the left and right sides of the equations below are equal:
\begin{align}
\nonumber [a_n,b_m]&=\delta_{n+m,0}, \\
 \label{4.2}[a_0,(\gamma b+\delta)^m]&=m\gamma (\gamma b+\delta)^{m-1}
\end{align}

We claim that the $h$-degrees of the vectors $a_{-n}$ and $b_{-n}$ are preserved under the adjoint action of $g$, namely
$(\pi(g)a)_{-n}$ is homogeneous and $h((\pi(g)a)_{-n})=h(a_{-n})$, similarly for $b_n$. 
Indeed, each term in the right side of (\ref{2.6}) has $h$-degree $-1$ which equals the $h$-degree of $a_{-n}$.
Similarly each term in the right side of (\ref{2.7}) has $h$-degree $1=h(b_{-n})$.

Suppose $ \gamma^{l} a_{-\lambda'} b_{-\mu'} (\gamma b+\delta)^{k}a_0^{s}$ is a term in $\pi(g) a_{-\lambda}b_{-\mu}\pi(g)^{-1}$, then 
\[
-p(\lambda)+p(\mu)=-l-p(\lambda')+p(\mu')-s.
\]
Hence $l=p(\lambda)-p(\mu)-p(\lambda')+p(\mu')-s=l(\lambda,\mu)-l(\lambda',\mu')-s$.

Moreover we can show that if $ \gamma^{l} a_{-\lambda'} b_{-\mu'} (\gamma b+\delta)^{k} a_0^s$ is a term in $\pi(g) a_{-\lambda}b_{-\mu}\pi(g^{-1})$, then $l+k$ must be $2p(\lambda)-2p(\mu)$, and $l \geq s$. 
This also can be proved by (\ref{2.6}) and (\ref{2.7}), and the fact that 
whenever moving the annihilation operators to the right, the sum of the indexes of $(\gamma b+\delta)$ and $\gamma$ will not be changed, and the index of $\gamma$ will not be decreased,
 thanks to the relation (\ref{4.2}).
Hence $k=2p(\lambda)-2p(\mu)-l=p(\lambda)-p(\mu)+p(\lambda')-p(\mu')+s=l(\lambda,\mu)+l(\lambda',\mu')+s$, and $s \leq [\frac{l(\lambda,\mu)-l(\lambda',\mu')}{2}]$.

Hence $c^s_{\lambda,\mu,\lambda',\mu'}$ is a constant.  And according to (\ref{2.6}) and (\ref{2.7}), powers of $\gamma$ and $(\gamma b+\delta)$ are the only information involved related to $g$ in the adjoint action, so the constant $c^s_{\lambda,\mu,\lambda',\mu'}$ has nothing to do with $g$. 

 The remaining properties are due to Lemma \ref{lemma2.3}.\qed

When (\ref{4.1}) is applied to $f(b)\in \mathcal O(\mathbb H)$, we get a refinement of (\ref{2.9}).

Now we will prove the remaining part of Theorem \ref{theorem2.5}, and we assume that $p(\lambda_0)=p(\mu_0)$.
 In this case it suffices to study the lifting of the constant function $ f(b)\equiv  1$, since the only modular form of weight $0$ is the constant function. 
We assume that there is a lifting of $1$ in $(V_{\lambda_0,\mu_0})^\Gamma_0$, and say 
\[
v=a_{-\lambda_0}b_{-\mu_0}+\sum_{\lambda,\mu: p(\lambda)<p(\mu)} a_{-\lambda}b_{-\mu}h_{\lambda,\mu}
\]
 is invariant under $g\in \Gamma$, for certain holomorphic functions $h_{\lambda,\mu} \in \mathcal O(\mathbb H)$.
Consider $\pi(g) v$, then the term corresponding to $a_{-\lambda}b_{-\mu}$ with $l(\lambda,\mu)=p(\lambda)-p(\mu)=-1$ only comes from $\pi(g)a_{-\lambda_0}b_{-\mu_0}$ and $\pi(g)a_{-\lambda}b_{-\mu}h_{\lambda,\mu}$  by Lemma \ref{lemma4.1}.
 So we have
 \begin{equation} \label{4.3}
(\gamma b+\delta)^{-2} h_{\lambda,\mu}(gb) -h_{\lambda,\mu}(b)+c_{\lambda_0,\mu_0,\lambda,\mu}^0 \gamma (\gamma b+\delta)^{-1}=0.
\end{equation}

Whenever $c^0_{\lambda_0,\mu_0,\lambda,\mu}$ vanishes, the function $h_{\lambda,\mu}$ is a modular form of weight $2$, 
and we may cancel the term $a_{-\lambda}b_{-\mu}h_{\lambda,\mu}$ by a lifting of $h_{\lambda,\mu}$ with the leading term $a_{-\lambda}b_{-\mu}h_{\lambda,\mu}$ as in Section 3. 
Hence we only need to consider the partition pairs with $c^0_{\lambda_0,\mu_0,\lambda,\mu}\neq 0$.

Recall that the Eisenstein series
\[
E_2(\tau):= 1+\frac{3}{\pi^2} \sum_{m\neq 0} \sum_{n\in\mathbb Z} \frac{1}{(m\tau +n)^2},
\]  
is holomorphic with the Fourier expansion
\[
E_2(\tau)=1-24\sum_{n=1}^\infty \sigma(n)q^n
\]
(cf. \cite{Z} p.19; \cite{A} p.69), where $q=e^{2\pi i \tau}$ and the divisor sum function $\sigma(n)=\sum_{d |n } d$.

 It is well-known that $E_2(\tau)$ is a quasi-modular form of weight $2$, with the transformation property (cf. \cite{Z} p.19; \cite{KZ})
\begin{equation}\label{4.4}
(\gamma\tau +\delta)^{-2}E_2(g\tau)=E_2(\tau) - \frac{6i}{\pi} \gamma (\gamma \tau +\delta)^{-1},\;\;\; \text{ for any } g\in SL(2,\mathbb Z).
\end{equation}

Obviously a rescaling of $E_2(b)$ by $\dfrac{\pi}{6i} c_{\lambda_0,\mu_0,\lambda,\mu}^0$, satisfies (\ref{4.3}), which is exactly the unique holomorphic (here the condition ``holomorphic" means holomorphic on $\mathbb H$ and all the cusps) solution up to modular forms.

Define $E(b):=\frac{\pi}{6i} E_2(b) $, then (\ref{4.4}) is equivalent to
\[
(\gamma\tau +\delta)^{-2}E(g\tau)=E(b) -  \gamma (\gamma \tau +\delta)^{-1}\;\;\; \text{ for any } g\in SL(2,\mathbb Z).
\]

Now we want to find an operator $A\in \bar{U}/K$, such that 
\[
\pi(g)  (a_{-\lambda_0} b_{-\mu_0} +AE(b)     )=a_{-\lambda_0}b_{-\mu_0} +AE(b),\;\;\;\;\text{ for any } g\in \Gamma.
\]

Notice that $\pi(g) AE(b)=(\pi(g)A\pi(g)^{-1}) E(gb)=(\pi(g)A\pi(g)^{-1})((\gamma b+\delta)^2E(b)-\gamma(\gamma  b+\delta))$,
it is natural to consider solutions of $A$, such that for $g\in SL(2,\mathbb R)$
\begin{align}
 \label{4.5}
\pi(g)A\pi(g)^{-1}&=A(\gamma b_0+\delta)^{-2},\\   
\label{4.6} A \gamma (\gamma b_0+\delta)^{-1} &=\pi(g) a_{-\lambda_0 }b_{-\mu_0}-a_{-\lambda_0} b_{-\mu_0}.
\end{align}

Similar to (\ref{3.2}), the equation (\ref{4.5}) can be replaced by the version of infinitesimal action, namely 
\begin{align}
\label{4.7}E_{(0)}. A&=0,\\
\label{4.8}H_{(0)}. A&=-2A,\\
\label{4.9}F_{(0)}. A&=2Ab.
\end{align}

\begin{lemma}
The condition (\ref{4.6}) implies (\ref{4.9}).
\end{lemma}
\noindent {\it Proof:} Notice that (\ref{4.9}) is equivalent to 
\[
F_{(0)}A b^i -A F_{(0)} b^i =2Ab^{i+1}, \;\;\;\text{ for all } i\geq 0,
\]
which is also equivalent to
\begin{equation}\label{4.10}
F_{(0)} A (-t b+1)^{-1} -A F_{(0)} (-t b+1)^{-1} = 2Ab(-tb+1)^{-1},
\end{equation}
where $t$ is a sufficiently small parameter.

Let $h(t):= \begin{pmatrix} 1 & 0 \\  -t & 1 \end{pmatrix}$. By (\ref{4.6}), we have
\[
A(-tb +1)^{-1} =-\frac{1}{t} (\pi(h(t)) a_{-\lambda_0}b_{-\mu_0} -a_{-\lambda_0}b_{-\mu_0}).
\]
Hence we have
\begin{align}\nonumber
F_0A(-tb+1)^{-1} & = - \left. \frac{1}{t} \frac{d}{ds} \right |_{s=0} \pi(h(s)) (\pi(h(t))a_{-\lambda_0}b_{-\mu_0}-a_{-\lambda_0}b_{-\mu_0}   )\\
\nonumber &=- \left. \frac{1}{t} \frac{d}{ds} \right |_{s=0}( (\pi(h(s+t)a_{-\lambda_0}b_{-\mu_0} -a_{-\lambda_0}b_{-\mu_0} )) -(\pi(h(s))a_{-\lambda_0}b_{-\mu_0} -a_{-\lambda_0}b_{-\mu_0} )     )\\
\nonumber &=- \left. \frac{1}{t} \frac{d}{ds} \right |_{s=0} (A(-s-t) ((-s-t) b+1)^{-1} -A(-s)(-sb+1)^{-1}    )\\
\label{4.11}&= A(-tb^2+2b)(-tb+1)^{-2}.
\end{align}

For any holomorphic function $f(b)\in \mathcal O(\mathbb H)$, 
\[
F_0 f(b)=\left.\frac{d}{ds} \right |_{s=0} \pi(h(s)) f(b) =\left. \frac{d}{ds} \right |_{s=0} f(\frac{b}{-sb+1})=f'(b) b^2.
\]
Hence we have
\begin{equation}\label{4.12}
A F_0(-tb+1)^{-1}= A tb^2(-tb+1)^{-2}.
\end{equation}

Using (\ref{4.11}) and (\ref{4.12}), the left hand side of (\ref{4.10}) equals
\[ 2Ab(-tb+1)^{-1},\]
which is exactly the right hand side.\qed

 Now we may assume 
 \begin{equation}\label{4.13}
 A=\sum_{(\lambda,\mu): p(\mu)>p(\lambda) } c_{\lambda,\mu} a_{-\lambda}b_{-\mu} a_0^{-l(\lambda,\mu)-1},
 \end{equation}
  where for any partition pair $(\lambda,\mu)$, $c_{\lambda,\mu}$ is a constant,  and $l(\lambda,\mu)$ as in (\ref{2.12}). Then $A$ satisfies (\ref{4.7}) and (\ref{4.8}). Applying (\ref{4.13}) to (\ref{4.6}), the left side becomes
 \[
 \sum_{(\lambda,\mu): p(\mu)>p(\lambda)} c_{\lambda,\mu}(-1)^{-l(\lambda,\mu)-1} (-l(\lambda,\mu)-1)! \gamma^{-l(\lambda,\mu)} a_{-\lambda}b_{-\mu} (\gamma b+\delta)^{l(\lambda,\mu)}.
 \]
which equals the right hand side of (\ref{4.6}), namely
\[
\sum_{(\lambda,\mu): p(\mu)>p(\lambda)} c_{\lambda_0,\mu_0,\lambda,\mu}^0 \gamma^{-l(\lambda,\mu)} a_{-\lambda}b_{-\mu}(\gamma b+\delta)^{l(\lambda,\mu)}.
\]
Comparing the coefficients, we have 
\[
c_{\lambda,\mu}=\dfrac{(-1)^{-l(\lambda,\mu)-1}}{(-l(\lambda,\mu)-1)!} c_{\lambda_0,\mu_0,\lambda,\mu}^0.
\]
Substituting the expression of $c_{\lambda,\mu}$ to (\ref{4.13}), we have 
\begin{equation}\label{4.14}
A=\sum_{(\lambda,\mu): p(\mu)>p(\lambda)} \dfrac{(-1)^{-l(\lambda,\mu)-1}}{(-l(\lambda,\mu)-1)!} c_{\lambda_0,\mu_0,\lambda,\mu}^0 a_{-\lambda} b_{-\mu} a_0^{-l(\lambda,\mu)-1},
\end{equation}
Thus, we obtain the following result.

\begin{theorem} \label{theorem4.3}
The operator
$A$ defined by (\ref{4.14}) satisfies (\ref{4.6})-(\ref{4.9}), hence $a_{-\lambda_0}b_{-\mu_0}+A E(b)$ is invariant under $\Gamma(1)=SL_2(\mathbb Z)$, and hence under arbitrary congruence subgroup $\Gamma$.
\end{theorem}

Theorem \ref{theorem4.3} completes the proof of Theorem \ref{theorem2.5}.
Apply Theorem \ref{theorem4.3}, we may calculate that the lifting of $1$ in $(V_{\lambda,\mu})^\Gamma_0$, for $\lambda=\mu=(1)$, equals
\[
 a_{-1}b_{-1}+2b_{-2}E(b)+b_{-1}^2 E'(b)=\omega+L_{-1} b_{-1}E(b),
 \]
 which is a Virasoro element in $\mathscr D^{\text{ch}}(\mathbb H,\Gamma)$ of central charge $2$. Therefore $\mathscr D^{\text{ch}}(\mathbb H,\Gamma)$ is a vertex operator algebra.

\section{Properties of vertex operator algebras $\mathscr D^{ch}(\mathbb H,\Gamma)$}

In this section, we study the structure of $\mathscr D^{\text{ch}}(\mathbb H,\Gamma)$ and calculate the character formula of $\mathscr D^{\text{ch}}(\mathbb H,\Gamma)$.

\begin{lemma} \label{lemma5.1}
Let $0\neq h \in M_k(\Gamma)$ for $k> 0$, so $h(z+N)=h(z)$ for some $N\in \mathbb Z _{>0}$. Suppose that $\sum_{i=0}^l g_i(z) h^{(i)}(z)$ is periodic with period $N$, where $g_i(z)\in \mathbb C[z]$ and $h^{(i)}(z)$ denotes $i$-th derivative of $h(z)$,
then $g_i(z)$ is a constant function for each $i$.
\end{lemma}

For simplicity, we denote by $D_f$ the space spanned by the lifting of modular forms $f\in M_{2k}(\Gamma), k\geq 1$ to $\mathscr D^{\text{ch}}(\mathbb H,\Gamma)$, namely 
\[
D_f:=\{  Af\in\mathscr D^{\text{ch}}(\mathbb H,\Gamma) \; |   \; A\in \bar{\mathcal U   }/K   \}.               
\]
And define
\[
D_1:=\text{Span}_\mathbb C \{  a_{-\lambda_0}b_{-\mu_0}+A E(b)\; |   \;  l(\lambda_0,\mu_0)=0 , A \text{ as in } (\ref{4.14}) \}.              
\]

The solution $A$ of (\ref{3.2}) (or equivalently (\ref{3.3})-(\ref{3.5})) 
 has the form 
\begin{equation}\label{5.1}
\sum_{\substack{ (\lambda,\mu)\leq (\lambda_0,\mu_0) \\  i\in \mathbb Z_{\geq 0}  }}    c_{\lambda,\mu}^i a_{-\lambda}b_{-\mu} a_0^i,
\end{equation}
where $c_{\lambda,\mu}^i\in \mathbb C$, and $c_{\lambda_0,\mu_0}^i$ vanishes unless $i=0$. 
Now we will show that for any $A\in \bar{\mathcal U}/K$, if $Af\in \mathscr D^{\text{ch}}(\mathbb H,\Gamma)$, where $f\in M_{2k}(\Gamma),k\geq 1$, then $A$ has the form (\ref{5.1}).
We assume that 
\[
A=\sum_{\substack{(\lambda,\mu)\leq (\lambda_0,\mu_0)  \\   i,j\in \mathbb Z_{\geq 0}            }    } c_{\lambda,\mu}^{i,j} a_{-\lambda}b_{-\mu} a_0^ib_0^j,
\]
 and our first step is to show that the expression of $A$ is free of $b_0$.
Indeed, since the congruence subgroup $\Gamma$ contains the translation matrix $\begin{pmatrix} 1 & N\\ 0 & 1 \end{pmatrix}$ for some positive integer $N$, $Af$ is fixed by $\begin{pmatrix} 1 & N\\ 0 & 1 \end{pmatrix}$,
which implies that 
\[
\sum_{i,j} c_{\lambda,\mu}^{i,j}a_0^i(b+N)^j f(b)=\sum_{i,j} c_{\lambda,\mu}^{i,j}a_0^i b^jf(b), \;\;\text{ for any partition pair } (\lambda,\mu). 
\] 
Fix $(\lambda,\mu)$ and let
$i_0$ be the maximal $i$, such that $ c^{i,j}_{\lambda,\mu} \neq 0 $, for some $ j$. Then $\sum_{j} c^{i_0,j}_{\lambda,\mu} b^j f^{(i_0)}(b)$ consists of the terms with highest derivation of $f(b)$, therefore 
according to Lemma \ref{lemma5.1}, all of $c^{i_0,j}_{\lambda,\mu}$ must be zero except possibly for $j=0$. 
By induction we can show that all of $c^{i,j}_{\lambda,\mu}=0$ except for $j=0$ and there would be no $b_0$ involved in the expression of $A$. Hence 
\begin{equation}\label{5.2}
A=a_{-\lambda_0}b_{-\mu_0} \sum_{i=0}^l c_{\lambda_0,\mu_0}^i a_0^i +\sum_{\substack{(\lambda,\mu)<(\lambda_0,\mu_0)\\ i\in \mathbb Z_{\geq 0}}} c_{\lambda,\mu}^i a_{-\lambda}b_{-\mu}a_0^i, \;\;\;\text{ with } c^l_{\lambda_0,\mu_0} \neq 0.
\end{equation}
 
Then we will show that $l=0$ in (\ref{5.2}). Comparing the terms corresponding to $a_{-\lambda_0}b_{-\mu_0}$ in both sides of $\pi(g)Af=Af$, we have 
\[
 \sum_{i=0}^l c_{\lambda_0,\mu_0}^i f^{(i)}(gb) (\gamma b+\delta)^{2l(\lambda_0,\mu_0)}= \sum_{i=0}^l c_{\lambda_0,\mu_0}^i f^{(i)}(b),
 \]
which implies $F(b):=\sum_{i=0}^l c_{\lambda_0,\mu_0}^i f^{(i)}(b)$ is a modular form of weight $-2l(\lambda_0,\mu_0)$. 
By induction we can show that 
\[
f^{(i)}(gb)= \sum_{j=0}^i c_{j}^i  \gamma ^{i-j} (\gamma b+\delta)^{2k+i+j} f^{(j)}(b), \;\;\;\text{ for } g\in \Gamma,
\]
where $c_j^i\in \mathbb Z_{>0}$ for $0 \leq j\leq i, \; c^i_i=1$. Hence, 
\[
\pi(g)F(b)=F(gb)=\sum_{i=0}^l\sum_{j=0}^i c_{\lambda_0,\mu_0}^ic_{j}^i  \gamma ^{i-j} (\gamma b+\delta)^{2k+i+j} f^{(j)}(b).
\]

Viewing $\pi(g)F(b)$ and $(\gamma b+\delta)^{-2l(\lambda_0,\mu_0)} F(b)$ as polynomials of $\gamma$, the leading terms of the two polynomials should be equal.
And the leading terms are 
\[
\sum_{j=0}^l c^l_{\lambda_0,\mu_0} c^l_j b^{2k+l+j}f^{(j)}(b) \gamma^{2k+2l},
\]
which is nonzero according to Lemma \ref{lemma5.1},
and $F(b)b^{-2l(\lambda_0,\mu_0)}\gamma ^{-2l(\lambda_0,\mu_0)}$ respectively.
Hence $l(\lambda_0,\mu_0)=-k-l$, and 
\[
\sum_{j=0}^l c^l_{\lambda_0,\mu_0} c^l_j b^{2k+l+j}f^{(j)}(b)= \sum_{i=0}^l c^{i}_{\lambda_0,\mu_0} b^{-2l(\lambda_0,\mu_0)} f^{(i)}(b),
\]
which is impossible unless
 $l=0$, and $-k=l(\lambda_0,\mu_0)$ by Lemma \ref{lemma5.1}.

\begin{lemma} \label{lemma5.2}
Let $A_if_i \in D_{f_i}$ be a nontrivial lifting of $f_i$ for $i=1,\cdots,n$, where $f_1,\cdots,f_n$ are linearly independent modular forms in $\oplus _{l\geq 1} M_{2l}(\Gamma)$. Then 
$A_1f_1,\cdots , A_nf_n$ are also linearly independent.
\end{lemma}
\noindent {\it Proof:} 
It suffices to prove the case for the modular forms of the same weight $2l$ with $l\geq 1$. Let $N>0$ be a fixed integer. 
Suppose $(\lambda_1,\mu_1) >\cdots>(\lambda_m,\mu_m)$ are all the partition pairs with the relation that $|\lambda_i|+|\mu_i|=N$ and $p(\lambda_i)-p(\mu_i)=-l$.

According to the discussion above Lemma \ref{lemma5.2} , we may assume $A_i$ has the leading term $a_{-\lambda_{k_i}}b_{-\mu_{k_i}}$,
and let
\[
A_i = \sum_{k_i\leq j \leq m }x_{ij} a_{-\lambda_j }b_{-\mu_j} +l.o.t,\;\;\;\text{ with } x_{ik_i}=1
\]
where we omit the lower order terms less than $a_{-\lambda_m}b_{-\mu_m}$ by $l.o.t$.
Notice that the term $a_{-\lambda_j}b_{-\mu_j} a_0^k $ for $j>k_i$ and $k>0$ may not appear in $A_i$, 
otherwise a combination of $A_i$ and the solutions $A$ in (\ref{3.2}) with $(\lambda_0,\mu_0)=(\lambda_{k_i},\mu_{k_i}),\cdots, (\lambda_{j+1},\mu_{j+1})$ respectively, will contradict the form (\ref{5.1}).

Suppose $\sum _{i=1}^n c_i A_if_i=0$, for some constants $c_i, i=1,\cdots,n$. Then 
\begin{equation}\label{equation5.2}
\sum_i c_i \sum_{ k_i\leq j\leq m} x_{ij} a_{-\lambda_j} b_{-\mu_j} f_i=0.
\end{equation}
 But the left hand side of (\ref{equation5.2}) equals $\sum_{j=1}^m \sum_{i:k_i\leq j} c_ix_{ij} a_{-\lambda_j}b_{-\mu_j} f_i=0$, hence $\sum_{i:k_i\leq j} c_ix_{ij} f_i=0$, for $1\leq j\leq m$.
 Because $f_1,\cdots,f_n$ are linear independent, $c_ix_{ij}=0$ for any $i$ such that $k_i\leq j$. So we have $c_i=0$, since $x_{ik_i}=1$. Therefore $c_i=0$ for $1\leq i\leq n$ and $A_1f_1,\cdots , A_nf_n$ are linearly independent. \qed

\begin{theorem}
Let $  \mathscr B$ be a homogeneous linear basis of modular forms in $\oplus_{k\geq 0} M_{2k}(\Gamma)$, then we have
\[
\mathscr D^{\text{ch}}(\mathbb H,\Gamma)=\oplus _{f\in \mathscr B} D_f.
\]
\end{theorem}

\

\ 
Now we will derive the character formula of $\mathscr D^{\text{ch}}(\mathbb H,\Gamma)$. 
The character formula of $\mathscr D^{\text{ch}}(\mathbb H,\Gamma)$ is the formal power series of variable $q$ defined by $\sum_{n=0}^\infty \dim \mathscr D^{\text{ch}}(\mathbb H,\Gamma)_n q^n$, which is 
$
tr\, q^{L_0}$.
We first consider that the trace $tr\, t^{H_{(0)}}q^{L_0}$ of the vertex subalgebra $\mathbb C[a_{-n},b_{-n} | n\geq 1]$. According to (\ref{3.8}), $tr\, t^{H_{(0)}}q^{L_0}$ is
\begin{equation}\label{5.3}
\Pi_{n=1}^{\infty} \dfrac{1}{1-t^2q^n} \Pi_{n=1}^\infty \dfrac{1}{1-t^{-2}q^n}=\sum_{n=0}^\infty \sum_{m=-\infty}^{\infty} c(m,n) q^nt^{m}
\end{equation} 
It is clear that $c(m,n)$ is the number of partition pairs $(\lambda,\mu)$ with $|\lambda|+|\mu|=n$ and $2(p(\lambda)-p(\mu))=m$.

As the character of $\mathscr D^{\text{ch}}(\mathbb H,\Gamma)$ coincides with the character of the graded algebra $\oplus_{\lambda,\mu} gr \mathscr D^{\text{ch}}(\mathbb H,\Gamma)_{\lambda,\mu}$, where $gr \mathscr D^{\text{ch}}(\mathbb H,\Gamma)_{\lambda,\mu}$ is the successive quotient $(V_{\lambda,\mu})^\Gamma_0/(V_{\lambda',\mu'})^\Gamma_0\cong M_{-2l(\lambda,\mu)}(\Gamma)$ by Theorem \ref{theorem2.5}, the character of $\mathscr D^{\text{ch}}(\mathbb H,\Gamma)$ equals
\begin{equation}\label{5.4}
\sum_{n=0}^\infty \sum_{m=-\infty}^\infty c(m,n) \dim M_{m}(\Gamma) q^n.
\end{equation}

\noindent  {\it Proof of Theorem \ref{theorem1.1}:}
When $m< 0$ or $m$ odd,  $\dim M_{m}(\Gamma)$ equals zero. So we let $m=2k$ for $k\geq 0$, and we will calculate $\sum_{n=0}^\infty c(2k,n) q^n$.

The left side of (\ref{5.3}) equals
\begin{gather*}
(1+\sum_{n,s\geq 1}  p_s(n)t^{2s}q^n )(1+\sum_{m,l \geq 1} p_l(m)t^{-2l}q^{m})\\
=1+\sum_{m,l\geq 1} p_l(m)(t^{-2l}+t^{2l})q^m+\sum_{m,n,l,s\geq 1}p_l(m)p_s(n)t^{-2l+2s}q^{m+n},
\end{gather*}
where $p_l(m)$ is the number of partitions of $m$ into exactly $l$ parts.

As $c(2k,n)$ is the coefficient of $q^nt^{2k}$ of the above formula, hence $\sum_{n=0}^\infty c(2k,n) q^n$ equals
\begin{equation}\label{5.5}
\sum_{n\geq 1}p_k(n)q^n +\sum_{l\geq 1} \sum_{m,n\geq 1} p_l(m)p_{l+k}(n) q^{m+n}
\end{equation}

Recall that a partition $\lambda$ has $k$ parts if and only if its conjugate partition $\lambda'$ has largest part $k$, 
where the conjugate partition $\lambda'$ is the partition whose Young diagram is obtained from interchanging rows and columns of $\lambda$.
So the generating function for partition with part $k$, is
\[
\sum_{n\geq 0}p_k(n)x^n=x^k\prod_{i=1}^k  \frac{1}{1-x^i}.
\]

Hence (\ref{5.5}) equals
\[
\pushQED{\qed}
\sum_{l\geq 0} q^{2l+k}\prod_{i=1}^l \dfrac{1}{1-q^i}\prod_{j=1}^{l+k} \dfrac{1}{1-q^j}.\qedhere
\popQED
\]

\section{Formulas for the lifting to a vertex subalgebra}

In this section, we consider the vertex algebra $\mathscr B:=\mathbb C[b_0,b_{-1},\cdots] \otimes _{\mathbb C[b]} \mathcal O(\mathbb H)$, generated by $ b_{-1}$ and $f(b)\in \mathcal O(\mathbb H)$.
 And we will find an explicit formula of a lifting of any nonconstant modular form of even weight to 
$\mathscr D^{\text{ch}}(\mathbb H,\Gamma(1)) \cap \mathscr B$.

As in Section 3, we prove that for any nonconstant modular form $f\in M_{2l}(\Gamma)$, 
there exists a lifting in $(V_{\lambda_0,\mu_0})^\Gamma_0$ for arbitrary partition pair $(\lambda_0,\mu_0)$ with the condition that $p(\mu_0)-p(\lambda_0)=l$.
Now we will give an explicit formula for a lifting of $f\in M_{2l}(\Gamma(1))$ to $(V_{\lambda_0,\mu_0})^{\Gamma(1)}_0\subset \mathscr D^{\text{ch}}(\mathbb H,\Gamma(1)) \cap \mathscr B$ 
for the case $\lambda_0=\emptyset$ and $p(\mu_0)=l$.
Since the action preserve the conformal weight, we will focus on $\mathscr B_n$ the weight $n$ subspace of $\mathscr B$.
Define a subspace $\mathscr B_n(l)\subset \mathscr B_n$ as follows
\[
\mathscr B_n(l) :=\text{Span}_\mathbb C\{  b_{-\mu} f(b) \in \mathscr B_n \;|\;   p(\mu)\geq l \}.
\]
Then we have a filtration which is stable under $\Gamma(1)$-action:
\[
\mathscr B_n=\mathscr B_n(0) \supset  \mathscr B_n(1) \supset \cdots \supset \mathscr B_n(n).
\]
Similar to the discussion in Section 2, $b_{-\mu}f(b)+\mathscr B_n(l+1)$, with $ p(\mu)=l$ is fixed by $\Gamma(1)$ if and only if
\begin{equation} \label{6.1}
f(gb)=(\gamma b_0+\delta)^{2l}f(b)
\end{equation}

Define $t^\mu_\nu:=(-1)^{p(\nu)-p(\mu)} c^0_{\emptyset, \mu,\emptyset, \nu}$.
Then (\ref{4.1}) turns into 
\begin{equation}\label{6.2}
\pi(g) b_{-\mu} f(b)=\sum_{\nu:p(\nu) \geq p( \mu)} t^{\mu}_{\nu}(-\gamma)^{-p(\mu)+p(\nu)} b_{-\nu} (\gamma b+\delta)^{-p(\mu)-p(\nu)}f(gb).
\end{equation}

Before we give the construction, we will explore some important properties of the integer coefficient $t^\mu_\nu$.

\begin{lemma} \label{lemma6.1}
Given any two partitions $\mu$ and $\nu$ of $n$ with the relation that $p(\mu) \leq p(\nu)$, we have
\begin{enumerate}
\item for $0\leq s\leq p(\nu)-p(\mu)$,
\[
 \sum_{\mu': p(\mu')=p(\mu)+s} t^\mu _{\mu'} t^{\mu'}_\nu= {p(\nu)-p(\mu) \choose  p(\mu') -p(\mu)} t^\mu_\nu.
 \]
\item for $s_i\geq 0$, with $i=1,2,\cdots,k-1$ and $\sum_{i=1}^{k-1} s_i \leq p(\nu)-p(\mu)$, 
\[
\sum_{\substack{\mu_1,\cdots,\mu_{k-1}:\\ p(\mu_i)=p(\mu_{i-1})+s_i}} t^{\mu}_{\mu_1}t^{\mu_1}_{\mu_2}\cdots t^{\mu_{k-1}}_{\nu}={p(\nu)-p(\mu) \choose p(\mu_1)-p(\mu),p(\mu_2)-p(\mu_1),\cdots,p(\mu_{k-1})-p(\mu_{k-2})}t^\mu_\nu,
\]
where $\mu_0$ is defined to be $\mu$.
\end{enumerate}
\end{lemma}

The above lemma is equivalent to the following one in a special form.
\begin{lemma} Given any two partitions $\mu$ and $\nu$ of $n$ with the relation that $p(\nu)=p(\mu)+k$, for $k\geq 2$, we have
\[
\sum_{\substack{\mu_1,\cdots,\mu_{k-1}: \\  p(\mu_i)=p(\mu)+i}  } t^\mu_{\mu_1} t^{\mu_1}_{\mu_2} \cdots t^{\mu_{k-1}}_\nu=k! t^\mu_\nu.
\]
\end{lemma}

\noindent {\emph Proof: }
Take $S_3$ to be an order three element in $PSL(2,\mathbb R)$, for example, let 
\[
S_3:=\begin{pmatrix} 
\frac{1}{2} & \frac{\sqrt 3}{2} \\
-\frac{\sqrt 3}{2}  & \frac{1}{2}
\end{pmatrix}
\]

We may compute the formula $\pi(S_3)^3 b_{-\mu}$ by iterating (\ref{6.2}), and it equals
\[
\sum_{\substack{\mu_1,\mu_2,\mu_3:\\ p(\mu)\leq p(\mu_1)\leq p(\mu_2) \leq p(\mu_3)}} 
(-1)^{-p(\mu_1)-p(\mu)} (\frac{\sqrt 3}{2})^{-p(\mu)+p(\mu_3)} t^\mu_{\mu_1} t^{\mu_1}_{\mu_2}t^{\mu_2}_{\mu_3}b_{-\mu_3} (-\frac{\sqrt 3}{2} b+\frac{1}{2})^{p(\mu_1)-p(\mu_3)}(-\frac{\sqrt 3}{2}b-\frac{1}{2})^{p(\mu)-p(\mu_2)} .
\]
On the other hand, $\pi(S_3)^3 b_{-\mu} =\pi(S_3^3) b_{-\mu}=b_{-\mu}$, since $S_3^3$ acts as identity operator.
So comparing the terms corresponding to $b_{-\nu}$ with $p(\nu)>  p(\mu)$ in the two expressions of $\pi(S_3)^3 b_{-\mu}$, we have for $b \in \mathbb H$ the following equation holds
\[
\sum_{\substack{\mu_1,\mu_2:\\ p(\mu)\leq p(\mu_1)\leq p(\mu_2) \leq p(\nu)}} 
(-1)^{p(\mu_1)-p(\mu)} t^\mu_{\mu_1} t^{\mu_1}_{\mu_2}t^{\mu_2}_{\nu} (-\frac{\sqrt 3}{2} b+\frac{1}{2})^{p(\mu_1)-p(\mu)}(-\frac{\sqrt 3}{2}b-\frac{1}{2})^{p(\nu)-p(\mu_2)} =0.
\]

Since the left side of the above equation is indeed a polynomial for $b$, it holds for arbitrary $b\in \mathbb C$. We take $b$ to be $0$, the equation becomes
\begin{equation} \label{6.3}
\sum_{\substack{\mu_1,\mu_2:\\ p(\mu)\leq p(\mu_1)\leq p(\mu_2) \leq p(\nu)}} 
 t^\mu_{\mu_1} t^{\mu_1}_{\mu_2}t^{\mu_2}_{\nu}(-\frac{1}{2})^{p(\mu_1)-p(\mu)+p(\nu)-p(\mu_2)} =0,
\end{equation}

Suppose the lemma holds for $k< l$, we will use (\ref{6.3}) to show the case $k=l$ by induction. Assume that $p(\nu)=p(\mu)+l$.
If $\{p( \mu_1),p(\mu_2)\}\subset \{p(\mu),p(\nu) \}$, then $t^\mu_{\mu_1} t^{\mu_1}_{\mu_2}t^{\mu_2}_{\nu}$
is exactly $t^\mu_\nu$ since whenever $p(\alpha)=p(\beta)$, then $t^{\alpha}_{\beta} =\delta_{\alpha,\beta}$, where $\delta_{\alpha,\beta}$ is the Kronecker delta function. 
If either $p( \mu_1)$ or $p(\mu_2)$ is distinct from both $p(\mu)$ and $p(\nu)$, we assume that $p(\mu_1)=p(\mu)+i$ and $p(\mu_2)=p(\mu_1)+j$.
Then running through all of $\mu_1$ and $\mu_2$ under the above condition,
the summation of $t^\mu_{\mu_1} t^{\mu_1}_{\mu_2}t^{\mu_2}_{\nu}(-\frac{1}{2})^{p(\mu_1)-p(\mu)+p(\nu)-p(\mu_2)} $ equals 
\[
 \frac{1}{i! j!(l-i-j)!} (-\frac{1}{2})^{  l-j}\sum_{\substack{\mu_1,\cdots,\mu_{l-1}: \\  p(\mu_s)=p(\mu)+s}  } t^\mu_{\mu_1} t^{\mu_1}_{\mu_2} \cdots t^{\mu_{l-1}}_\nu,
\]
where we use the induction assumption.

Thanks to the combinatory equation
\[
\sum_{i=0}^l \sum_{j=0}^{l-i} \frac{1}{i! j!(l-i-j)!} (-\frac{1}{2})^{  l-j} =0,
\]
 the left side of (\ref{6.3}) equals
\[
(2(-\frac{1}{2})^l+1)t^\mu_\nu-\frac{1}{l!}(2(-\frac{1}{2})^l+1)\sum_{\substack{\mu_1,\cdots,\mu_{l-1}: \\  p(\mu_s)=p(\mu)+s}  } t^\mu_{\mu_1} t^{\mu_1}_{\mu_2} \cdots t^{\mu_{l-1}}_\nu
\]
hence we have proved the case $k=l$.\qed

Now suppose $f$ is a modular form of weight $2l$, we will derive a formula for the lifting of $f$ in $\mathscr D^{\text{ch}}(\mathbb H,\Gamma) \cap \mathscr B_n(l)$ with the leading term $b_{-\mu}f(b)$,
where $p(\mu)=l$, and $|\mu|=n$.

We define 
\begin{equation}\label{6.4}
F_k(b)=\sum_{s=0}^k \sum_{\mu': p(\mu')=p(\mu)+s} c_{\mu'} b_{-\mu'} f^{(s)}(b),
\end{equation}
where $f^{(s)}(b)$ denotes $s$-th derivatives of $f(b)$, $c_{\mu'}:= t^\mu_{\mu'}  \cdot  \Pi_{t=0}^{p(\mu')-l-1}  (2l+t) ^{-1 } \in \mathbb R$,
 for partition $\mu'$ such that $p(\mu') >l$ and $c_\mu=1$. Oviously, $F_k(b)\in \mathscr B_n(l)$, for $k\geq 0$, and $F_0(b)=b_{-\mu} f(b)$.

Hence we have for any $g=\begin{pmatrix} \alpha & \beta\\ \gamma & \delta \end{pmatrix}\in \Gamma(1)$, obviously 
$\pi(g) F_0(b)-F_0(b) \in \mathscr B_{n}(l+1).$
Suppose that $\pi(g) F_{k-1}(b)-F_{k-1}(b)\in \mathscr B_n(l+k)$, we will show that 
\begin{equation}\label{6.4}
\pi(g) F_{k}(b)-F_{k}(b)\in \mathscr B_n(l+k+1).
\end{equation}
Note that $\pi(g) F_{k-1}(b)-F_{k-1}(b)$ equals
\begin{equation}\label{6.6}
\sum_{s=0}^{k-1} \sum_{\substack{\mu',\nu:\\ p(\mu')=p(\mu)+s}}c_{\mu'} t^{\mu'}_\nu (-\gamma)^{p(\nu)-p(\mu')} (\gamma b+\delta)^{-(p(\nu)+p(\mu'))} b_{-\nu} f^{(s)}(gb)-\sum_{s=0}^ {k-1} \sum_{\mu':p(\mu')=p(\mu)+s} c_{\mu'} b_{-\mu'} f^{(s)}(b)
\end{equation}

As (\ref{6.6}) is contained in $\mathscr B_n(l+k)$, the term $b_{-\nu} h(b)$, $0\neq h(b)\in\mathcal O(\mathbb H)$ will not appear for any  partition $\nu$ such that $p(\nu)=l+k-1$.
Therefore we have
\[
\sum_{s=0}^{k-1} \sum_{\mu':p(\mu')=l+s} c_{\mu'} t^{\mu'}_\nu (-\gamma)^{p(\nu)-p(\mu')}(\gamma b+\delta)^{-p(\nu)-p(\mu')} f^{(s)}(gb)-c_\nu f^{(k-1)}(b)=0
\]
Taking derivative towards $b$ leads to the equation that
\begin{align*}
&\sum_{s=0}^{k-1} \sum_{\mu':p(\mu')=l+s} c_{\mu'} t^{\mu'}_\nu(2l+k+s-1) (-\gamma)^{k-s}(\gamma b+\delta)^{-2l-k-s} f^{(s)}(gb)\\
&+\sum_{s=0}^{k-1} \sum_{\mu': p(\mu')=l+s} c_{\mu'} t^{\mu'}_\nu (-\gamma)^{k-s-1} (\gamma b+\delta)^{-2l-k-s-1} f^{(s+1)}(gb)-c_\nu f^{(k)}(b)=0
\end{align*}
We substitute $s$ by $s-1$ in the second summation above, the equation becomes
\begin{align*}
& (2l+k-1) t^{\mu}_\nu (-\gamma)^k (\gamma b+\delta)^{-2l-k} f(gb)\\
& +\sum_{s=1}^{k-1}\left(\sum_{\mu' :p(\mu')=l+s} (2l+k+s-1)c_{\mu'} t^{\mu'}_\nu  + \sum_{\mu': p(\mu')=l+s-1} c_{\mu'} t^{\mu'}_\nu   \right) (-\gamma )^{k-s} (\gamma b+\delta)^{-2l-k-s} f^{(s)}(gb)\\
& +c_\nu (\gamma b+\delta)^{-2l-2k} f^{(k)}(gb) -c_\nu f^{(k)}(b)=0
\end{align*}
Fix an arbitrary partition $\chi$ of $n$ with $p(\chi)=l+k$. Then multiply $t^\nu_\chi$ to the above equation and take sum of all partition $\nu$ of $n$ with $p(\nu)=l+k-1$. We have the equation below

\begin{align*}
&\sum_{\nu :p(\nu)=l+k-1} (2l+k-1)  t^\nu_\chi t^\mu_\nu(-\gamma)^k (\gamma b+\delta)^{-2l-k} f(gb)\\
&+ \sum_{s=1}^{k-1} \left(\sum_{\substack{\mu',\nu :\\ p(\mu')=l+s\\p(\nu)=l+k-1}} (2l+k+s-1)c_{\mu'}  t^{\mu'}_\nu t^\nu_\chi +\sum_{\substack{\mu',\nu:\\p(\mu')=l+s-1\\ p(\nu) =l+k-1}} c_{\mu'} t^{\mu'}_\nu t^\nu_\chi  \right)(-\gamma)^{k-s} (\gamma b+\delta)^{-2l-k-s} f^{(s)}(b)\\
&+ \sum_{\nu:p(\nu)=l+k-1} c_\nu t^\nu_\chi (\gamma b+\delta )^{-2l-2k} f^{(k)}(gb) - \sum_{\nu: p(\nu)=l+k-1}c_\nu t^\nu_\chi f^{(k)}(b)=0
\end{align*}

The left side of the above equation, according to Lemma \ref{lemma6.1} and definition of $c_\mu$, equals
\begin{align*}
& (2l+k-1) k t^\mu_\chi (-\gamma )^k (\gamma b+\delta)^{-2l-k} f(gb)\\
&+\sum_{s=1}^{k-1} {k \choose s}  k(2l+k-1)\Pi_{t=0}^{s-1} (2l+t)^{-1} t^\mu_{\chi} (-\gamma)^{k-s} (\gamma b+\delta)^{-2l-k-s} f^{(s)}(gb)\\
&+ k\Pi_{t=0}^{k-2} (2l+t)^{-1} t^\mu _{\chi} (\gamma b_0 +\delta)^{-2l-k} f^{(k)}(gb)-k\Pi_{t=0}^{k-2}(2l+t)^{-1} t^\mu_\chi f^{(k)}(b)
\end{align*}
which is exactly the function corresponding to $b_{-\chi}$ in  $g F_k(b)-F_k(b)$. 
Hence $gF_{n-l}(b)-F_{n-l}(b) \in \mathscr B_n(n+1)=0$, thus $F_{n-l}(b)$ defined in (\ref{6.4}) is a lifting of $f$ to $\mathscr D^{\text{ch}}(\mathbb H,\Gamma) \cap \mathscr B^{\Gamma(1)}_0$ with the leading term $b_{-\mu}f(b)$.

\

\

\end{document}